\documentclass{amsart}

\usepackage{m-pictex, epsfig, amsthm, amssymb, amsmath, amsxtra, graphicx, psfrag, latexsym}
\usepackage[all]{xy}
\usepackage{enumerate}

\DeclareMathOperator{\fix}{Fix}
\DeclareMathOperator{\inti}{Int}
\DeclareMathOperator{\sta}{Stab}

\DeclareMathOperator{\gt}{\Gamma_3}
\DeclareMathOperator{\g}{\Gamma}

\DeclareMathOperator{\iso}{Isom}

\DeclareMathOperator{\f}{\mathcal F}
\DeclareMathOperator{\vc}{\mathcal V\mathcal C}
\DeclareMathOperator{\vcm}{\mathcal V\mathcal C_{\infty}}
\DeclareMathOperator{\fin}{\mathcal F\mathcal I\mathcal N}
\DeclareMathOperator{\all}{\mathcal A\mathcal L\mathcal L}
\DeclareMathOperator{\tr}{\mathcal T\mathcal R}

\theoremstyle{plain}
\newtheorem{theorem}{Theorem}[section]
\newtheorem{lemma}[theorem]{Lemma}

\newtheorem{proposition}[theorem]{Proposition}
\newtheorem{remark}[theorem]{Remark}

\newtheorem{claim}[theorem]{Claim}

\theoremstyle{definition}
\newtheorem{definition}[theorem]{Definition}

\begin{document}

\title[Classifying spaces and lower algebraic K-theory.]
      {Relative hyperbolicity, classifying spaces, and lower algebraic K-theory.}
\author{Jean-Fran\c{c}ois\ Lafont}
\address{Department of Mathematics\\
         Ohio State University\\
         Columbus, OH  43210}
\email[Jean-Fran\c{c}ois\ Lafont]{jlafont@math.ohio-state.edu}
\author{Ivonne J.\ Ortiz}
\address{Department of Mathematics and Statistics\\
         Miami University\\
         Oxford, OH 45056}
\email[Ivonne J.\ Ortiz]{ortizi@muohio.edu} 
%\subjclass[2000]{19B28,19A31, 19D35} 
%\keywords{Relative hyperbolic groups, classifying
%spaces for a family of subgroups, lower algebraic $K$-theory}
%\date{\today}

\begin{abstract}
For $\Gamma$ a relatively hyperbolic group, we construct a model for
the universal space among $\Gamma$-spaces with isotropy on the
family $\vc$ of virtually cyclic subgroups of $\Gamma$.  We provide
a recipe for identifying the maximal infinite virtually cyclic
subgroups of Coxeter groups which are lattices in $O^+(n,1)=
\iso(\mathbb H^n)$. We use the information we obtain to explicitly
compute the lower algebraic K-theory of the Coxeter group $\gt$
(a non-uniform lattice in $O^+(3,1)$). Part of this computation
involves calculating certain Waldhausen Nil-groups for $\mathbb Z[D_2]$, $\mathbb Z [D_3]$.
\end{abstract}

\maketitle
\section[Introduction]{Introduction}

Let $\g$ be a discrete group and  let $\f$ be a family  of subgroups
$\g$. A $\Gamma$-CW-complex $E$ is a model for the classifying space
$E_{\f}(\g)$ with isotropy on $\f$ if the $H$-fixed point sets $E^H$
are contractible for all $H \in \mathcal F$ and empty otherwise.  It
is characterized by the universal property that for every
$\Gamma$-CW complex $X$ whose isotropy groups are all in $\mathcal
F$, one can find an equivariant continuous map $X \rightarrow
E_{\f}(\g)$ which is unique up to equivariant homotopy.  The two
extreme cases are $\mathcal F=\all$, where $E_{\f}(\g)$ can be taken
to be a point, and $\mathcal F=\tr$, where $E_{\f}(\g)$ is a model
for $E\Gamma$.

For the family of finite subgroups, the space $E_{\f}(\g)$ has  nice geometric models for 
various classes of groups $\g$.  For instance, in the case where $\g$ is is a
discrete subgroup of a virtually connected Lie group \cite{FJ93}, where $\g$ is
word hyperbolic group \cite{MS02}, an arithmetic group \cite{BS73} \cite{S79}, 
the outer automorphism group of a free group \cite{CV86}, a  mapping class groups \cite{K83}, 
or a one relator group \cite{LS77}.  For a thorough
survey on classifying spaces, we refer the reader to L\"uck \cite{Lu04}.

One motivation for the study of these classifying spaces comes from
the fact that they appear in the Farrell-Jones Isomorphism
Conjecture about the algebraic $K$-theory of groups rings (see  \cite{FJ93}). Because
of this conjecture the computations of the relevant $K$-groups can
be reduced to the computation of certain equivariant homology groups
applied to these classifying spaces for the family of finite groups
$\fin$ and the family of virtually cyclic subgroups $\vc$ (where a
group is called {\it virtually cyclic} if it has a cyclic subgroup
of finite index).  In this paper we are interested in classifying
spaces with isotropy in the family $\vc$ of virtually cyclic
subgroups.

We start out by defining the notion of an adapted family associated
to a pair of families $\mathcal F \subset \mathcal G$ of subgroups.
We then explain how, in the presence of an adapted family, a classifying
space $E_{\mathcal F} (\Gamma)$ can be modified to obtain a classifying
space $E_{\mathcal G} (\Gamma)$ for the larger family.
In the situation we are interested in, the smaller family will be $\fin$ and
the larger family will be $\vc$.

Of course, our construction is only of interest if we can find
examples of groups where there is already a good model for
$E_{\fin}(\g)$, and where an adapted family can easily be found. For
$\g$ a relatively hyperbolic group in the sense of Bowditch
\cite{Bo} (or equivalently relatively hyperbolic with the bounded
coset penetration property in the sense of Farb \cite{Fa98}),
Dahmani has constructed a model for $E_{\fin}(\g)$.  We show that
the family consisting of all conjugates of peripheral subgroups,
along with all maximal infinite virtually cyclic subgroups {\it not}
conjugate into a peripheral subgroup, forms an adapted family for
the pair $(\fin, \vc)$. Both the general construction, and the
specific case of relatively hyperbolic groups, are discussed in
Section 2 of this paper.

In order to carry out our construction of the classifying spaces for the family
$\vc$ for these groups, we need to be able to classify the maximal
infinite virtually cyclic subgroups.  We establish a systematic
procedure to complete this classification for arbitrary Coxeter
groups arising as lattices in $SO(n,1)$. We next focus on the group
$\gt$, a Coxeter group which is known to be a non-uniform lattice in
$SO(3,1)$.  In this specific situation, it is well known that the
action of $\gt$ on $\mathbb H^3$ is a model for $E_{\fin}(\gt)$, and
that the group $\gt$ is hyperbolic relative to the cusp group (in
this case the 2-dimensional crystallographic group $P4m$).  Our
construction now yields an 8-dimensional classifying space for
$E_{\vc}(\gt)$.  These results can be found in Section 3 of our
paper.

Since the Farrell-Jones isomorphism conjecture is known to hold for
lattices in $SO(n,1)$, we can use our 8-dimensional classifying
space for $E_{\vc}(\gt)$ to compute the lower algebraic K-theory of
(the integral group ring of) $\gt$.  The computations are carried
out in Section 4 of the paper, and yields an explicit result for
$K_n(\mathbb Z \gt)$ when $n\leq -1$.  The $\tilde K_0(\mathbb Z \gt)$
and $Wh(\gt)$ terms we obtain involve some Waldhausen Nil-groups.

In general, very little is known about Waldhausen Nil-groups. In
Section 5, we provide a complete \emph{explicit} determination of the
Waldhausen Nil-groups that occur in $\tilde K_0(\mathbb Z \gt)$ and
$Wh(\gt)$.  The approach we take was suggested to us by F.T.
Farrell, and combined with the computations in Section 4, yields the
first example of a lattice in a semi-simple Lie group for which (1) the lower 
algebraic K-theory
is \emph{explicitly} computed, but (2) the relative assembly map induced by
the inclusion $\fin\subset \vc$ is {\it not} an isomorphism.  The
result of our computations can be summarized in the following:

\begin{theorem}

Let $\Gamma_{3}=O^{+}(3,1) \cap GL(4,\mathbb Z)$. Then the lower
algebraic $K$-theory of the integral group ring of $\gt$ is
given as follows:
\[
\begin{aligned}
Wh(\gt) &\cong \bigoplus _\infty \mathbb Z /2\\
\tilde{K}_0(\mathbb Z \gt) &\cong
\mathbb Z/4 \oplus \mathbb Z/4 \oplus \bigoplus _\infty \mathbb Z /2\\
K_{-1}(\mathbb Z \gt) &\cong \mathbb Z \oplus \mathbb Z, \; and\\
K_n(\mathbb Z \gt) &\cong 0, \; for \; n < -1.
\end{aligned}
\]
where the expression $\bigoplus _\infty \mathbb Z /2$ refers to a
countable infinite sum of $\mathbb Z/2$.
\end{theorem}

Finally, we note that most of the results in this paper apply in a
quite general setting, and in particular to any Coxeter group that
occurs as a lattice in $O^+(n,1)$.  Due to space constraints, we
have only included the computations for the group $\gt$.  In a
forthcoming paper, the authors will carry out the corresponding
computations for the lower algebraic K-theory of the remaining
3-simplex hyperbolic reflection groups.

\vskip 10pt

\section[The classifying space]{A Model for $E_{\vc}(\g)$}

Let $\Gamma$ be a discrete group and $\mathcal F$ be a  family of
subgroups of $\Gamma$ closed under inclusion and conjugation, i.e.\ if
$H \in \mathcal F$ then $gH'g^{-1} \in \mathcal F$ for all $H' \subset H$
and all $g \in \Gamma$. Some examples for $\mathcal F$ are $\mathcal
T\mathcal R$, $\mathcal F\mathcal I\mathcal N$, $\mathcal V\mathcal C
$, and $\mathcal A\mathcal L\mathcal L$, which are the families
consisting of the trivial group, finite subgroups, virtually cyclic
groups, and all subgroups respectively.

\begin{definition}
Let $\Gamma$ be any finitely generated group, and $\mathcal F \subset
\widetilde{\mathcal F}$ a pair of families of subgroups of $\g$, we say that
a collection $\{H_{\alpha}\}_{\alpha \in I}$ of subgroups of $\g$ is
{\it adapted} to the pair $(\mathcal F, \widetilde{\mathcal F})$ provided
that:
\begin{enumerate}
\item For all $G, H \in \{H_{\alpha}\}_{\alpha \in I}$,
either $G=H$, or $G \cap H \in \mathcal F$.

\vspace{.2cm}
\item The collection $\{H_{\alpha}\}_{\alpha \in I}$ is {\it conjugacy
closed} i.e.\ if $G \in \{H_{\alpha}\}_{\alpha \in I}$ then
$gGg^{-1} \in \{H_{\alpha}\}_{\alpha \in I}$ for all $g \in \g$.

\vspace{.2cm}
\item Every $G \in \{H_{\alpha}\}_{\alpha \in I}$ is {\it self-normalizing}, i.e.\ $N_{\g}(G)=G$.

\vspace{.2cm}
\item For all $G \in \widetilde{\mathcal F} \setminus \mathcal F$, there
exists $H \in \{H_{\alpha}\}_{\alpha \in I}$ such that $G \leq H$.
\end{enumerate}
\end{definition}

\begin{remark}
The collection $\{\g\}$ consisting of just $\g$ itself is adapted to
every pair $(\mathcal F, \widetilde{\mathcal F})$ of families of
subgroups of $\g$. Our goal in this section is to show how, starting
with a model for $E_{\f}(\g)$, and a collection
$\{H_{\alpha}\}_{\alpha \in I}$ of subgroups adapted to the pair
$(\mathcal F, \widetilde{\mathcal F})$ , one can build a model for
$E_{\widetilde{\mathcal F}}(\g)$.
\end{remark}

\subsection[the construction]{The Construction}

\begin{enumerate}
\item For each subgroup $H$ of  $\Gamma$, define the induced family of
subgroups $\widetilde{\mathcal F}_H$ of $H$ to be
$\widetilde{\mathcal F}_H :=\{F \cap H|F \in \widetilde{\mathcal
F}\}$.  Note that if $g \in \g$, conjugation by $g$ maps $H$ to
$g^{-1}Hg \leq \g$, and sends $\widetilde{\f}_H$ to
$\widetilde{\f}_{g^{-1}Hg}$.

Let $E_H$ be a model for the classifying space
$E_{\widetilde{\mathcal F}_H}(H)$ of $H$ with isotropy in
$\widetilde{\f}_H$. Define a new space $E_{H;\g} = \coprod_{\g/H}
E_H$. This space consists of the disjoint union of copies of $E_H$,
with one copy for each left-coset of $H$ in $\g$. Note that $E_H$ is
contractible, but $E_{H;\g}$ is not (since it is not
path-connected).

\item Next, we define a $\g$-action on the space $E_{H;\g}$. Observe
that each component of $E_{H;\g}$ has a natural $H$-action; we want
to ``promote" this action to a $\g$-action. By abuse of notation,
let us denote by $E_{gH}$ the component of $E_{H;\g}$ corresponding
to the coset $gH \in \g/H$. Fix a collection $\{g_i H | i \in I\}$
of left-coset representatives, so that we now have an identification
$E_{H;\g} = \coprod_{i \in I} E_{g_iH}$. Now for $g \in \g$, we
define the $g$-action on $E_{H;\g}$ as follows: g maps each
$E_{g_iH}$ to $E_{{gg_i}H} \cong E_{g_jH}$, for some $j \in I$.
Recall that both $E_{g_iH}$ and $E_{g_jH}$ are copies of $E_H$, and
that $g^{-1}_j g g_i \in H$; since $H$ acts on $E_H$, we define the
$g$-action from $E_{g_iH}$ to $E_{g_jH}$ to be:

\[
\xymatrix{\ar @{} [dr] |{\circlearrowleft}
 E_{g_iH} \ar[d]_{\cong} \ar[r]^g & E_{g_jH}
                 \ar[d]^{\cong} \\
E_H \ar[r]_{g^{-1}_j g g_i} & E_H}
\]

\noindent that is the $g^{-1}_j g g_i$-action via the identification
$E_{g_iH} \cong E_H \cong E_{g_jH}$. Note that while the $\g$-action
one gets out if this depends on the choice of left-coset
representatives, different choices of representatives yield
conjugate actions, i.e.\ there exist a $\Gamma$-equivariant
homomorphism from  $\coprod_{i \in I} E_{g_iH}$ and  $\coprod_{i \in
I} E_{g'_iH}$.

\item Now assume that we have a pair of families $\mathcal F \subseteq
\widetilde{\mathcal F}$, and a collection  $\{H_{\alpha}\}_{\alpha
\in I}$ of subgroups of $\g$ adapted to the pair $(\mathcal F,
\widetilde{\mathcal F})$. Select a subcollection
$\{\Lambda_{\alpha}\}_{\alpha \in \widetilde{I}} \subseteq
\{H_{\alpha}\}_{\alpha \in I}$ consisting of one subgroup from each
conjugacy class in $\{H_{\alpha}\}_{\alpha \in I}$. We form the
space $\widehat{X}$ by taking the join of $X$, a model for
$E_{\mathcal F}(\g)$, with the disjoint union of the spaces $E_{\Lambda_{\alpha}; \g}$,
defined as in (1), with the naturally induced $\g$-action given in
(2), i.e.\ $\widehat X = X \ast \coprod_{\alpha \in \widetilde{I}}
E_{\Lambda_{\alpha}; \g}$.

\end{enumerate}

\begin{theorem}
The space $\widehat{X}$ is a model for $E_{\widetilde{\mathcal F}}(\g)$.
\end{theorem}

\begin{proof}
We start by noting that if a space $A$ is contractible, and $B$ is
any space, then $A \ast B$ is contractible; in particular the space
$\widehat X$ is contractible (since it is a join with the
contractible space X).

We need to show two points:

\begin{enumerate}
\item $\widehat{X}^H$ is contractible if $H \in \widetilde{\f}$.

\item $\widehat{X}^H=\emptyset$ if $H \notin \widetilde{\f}$.
\end{enumerate}

Let us concentrate on the first point; assume $H \in \widetilde{
\f}$. Note that if $H \in \f$, then since $X$ is a model for
$E_{\f}(\g)$, we have that $X^H$ is contractible, and since
$\widehat{X}^H=(X \ast \coprod_{\alpha \in
\widetilde{I}}E_{\Lambda_{\alpha}; \g})^H=X^H \ast \coprod_{\alpha \in
\widetilde{I}}E^H_{\Lambda_{\alpha}; \g}$, we conclude that
$\widehat{X}^H$ is contractible.

Now assume $H \in \widetilde{\f} \setminus \f$. From property (4) of
an adapted family for the pair $(\widetilde{\f}, \f)$ (see
Definition 2.1), there exists $\Lambda_{\alpha}$ for some $\alpha
\in \widetilde{I}$, such that $H$ can be conjugated into
$\Lambda_{\alpha}$. We now make the following claim:

\begin{claim}
$E^H_{\Lambda_{\alpha}; \g}$ is contractible, and
$E^H_{\Lambda_{\beta}; \g}=\emptyset$, for all $\beta \neq \alpha$.
\end{claim}

Note that if we assumed Claim 2.4, we immediately get that
$\widehat{X}^H$ is contractible. Indeed, since $X$ is a model for
$E_{\f}(\g)$ and $H \notin \f$, $X^H=\emptyset$, and we obtain:
\[
\widehat{X}^H =E^H_{\Lambda_{\alpha}; \g} \ast X^H \ast \coprod_{\beta \neq \alpha}
E^H_{\Lambda_{\beta}; \g}=E^H_{\Lambda_{\alpha}; \g} \ast \emptyset
\ast \coprod_{\beta \neq \alpha} \emptyset \cong E^H_{\Lambda_{\alpha}; \g}.
\]

\begin{proof}[Proof of Claim 2.4]
We first show that, in $E_{\Lambda_{\alpha};\g} = \coprod_{i \in I}
E_{g_i\Lambda_{\alpha}}$, there exists a connected component that is
fixed by $H$. Since $H$ is conjugate to $\Lambda_{\alpha}$ by some
$k \in \g$, consider the $\Lambda_{\alpha}$ coset
$g_i\Lambda_{\alpha}$ containing $k$. We claim that for all $h \in
H$, we have $h \cdot g_i\Lambda_{\alpha}=g_i\Lambda_{\alpha}$.
First observe the following easy criterion for an element $h\in \Gamma$
to map a coset $g_i\Lambda_\alpha$ to itself:
\[
h \cdot g_i\Lambda_{\alpha}=g_i\Lambda_{\alpha} \Longleftrightarrow
g^{-1}_i h g_i \cdot \Lambda_{\alpha}=\Lambda_{\alpha}
\Longleftrightarrow g^{-1}_i h g_i \in \Lambda_{\alpha}.
\]
Now in our specific situation, by the choice of $g_i\Lambda_{\alpha}$, we have
$k=g_i\widetilde{k} \in g_i\Lambda_{\alpha}$, where $k$ conjugates
$h$ into $\Lambda_{\alpha}$, and $\widetilde{k} \in
\Lambda_{\alpha}$. Since $g_i=k\widetilde{k}^{-1}$, substituting we
see that $g^{-1}_i h g_i=\widetilde{k}(k^{-1} h k)\widetilde{k}^{-1}$.  But
by construction, we have that $k^{-1} h k\in \Lambda_{\alpha}$, and since
$\widetilde{k}\in \Lambda_{\alpha}$, we conclude that $g^{-1}_i h g_i
\in \Lambda_{\alpha}$. So by the criterion above, we see that indeed, every element of
$H$ maps the component  $E_{g_i\Lambda_{\alpha}}$ to itself.

Our next step is to note that an element $g \in H$ with $g \neq 1$
can map {\it at most one} of the $E_{g_i\Lambda_{\alpha}}$ to
itself. In fact, if $g \in H$ maps $E_{g_i\Lambda_{\alpha}}$ to
itself and $E_{g_j\Lambda_{\alpha}}$ to itself, then from the
definition of the the action given in step (2) of the construction,
we have that $H \subset g_i \Lambda_{\alpha} g^{-1}_i$, and $H
\subset g_j \Lambda_{\alpha} g^{-1}_j$. From condition (1) in
Definition 2.1, we have that either $g_i \Lambda_{\alpha} g^{-1}_i =
g_j \Lambda_{\alpha} g^{-1}_j$, or $H \subset g_i \Lambda_{\alpha}
g^{-1}_i \cap g_j \Lambda_{\alpha} g^{-1}_j\in \f$. But this second
possibility can not occur, as $H \in \widetilde{\f} \setminus \f$.
Hence we have that $g_i \Lambda_{\alpha} g^{-1}_i = g_j
\Lambda_{\alpha} g^{-1}_j$. Therefore $(g^{-1}_j g_i)
\Lambda_{\alpha} (g^{-1}_i g_j)=\Lambda_{\alpha}$, that is,
$g^{-1}_j g_i \in N_{\g}(\Lambda_{\alpha})$. But the
self-normalizing condition (see condition (3) in Definition 2.1)
forces $g^{-1}_j g_i \in \Lambda_{\alpha}$, thus $g_i \in
g_j\Lambda_{\alpha}$, i.e. the cosets $g_i\Lambda_{\alpha}$, and
$g_i\Lambda_{\alpha}$ had to be the same to begin with.

Up to this point, we know that the $H$-action on
$E_{\Lambda_{\alpha}; \g}$ maps the $E_{g_i \Lambda_{\alpha}}$
component to itself, and permutes all the others. This implies
$E^H_{\Lambda_{\alpha};\g} =E^H_{g_i \Lambda_{\alpha}}$. But the
$H$-action on $E_{g_i \Lambda_{\alpha}}$ coincides with the
$g^{-1}_i H g_i$-action on $E_{\Lambda_{\alpha}}$ (see step (2) of
the construction). From step (1) of the construction, since
$g^{-1}_i H g_i \in \widetilde{\f}_{\Lambda_{\alpha}}$, and
$E_{\Lambda_{\alpha}}$ is a model for
$E_{\widetilde{\f}_{\Lambda_{\alpha}}}(\Lambda_{\alpha})$, we have
that $E^H_{g_i \Lambda_{\alpha}}=E^H_{ \Lambda_{\alpha}; \g}$ is
contractible.
\end{proof}

This completes the verification of the first point, i.e.\ $\widehat
X^H$ is contractible if $H \in \widetilde{\f}$.

We now verify the second point, i.e.\ $\widehat{X}^H=\emptyset$ if
$H \notin \widetilde{\f}$. But this is considerably easier. In fact,
if $H \notin \widetilde{\f}$, then $H \notin \f$, and
$X^H=\emptyset$; so let us focus on the $H$-action on the individual
$E_{\Lambda_{\alpha};\g}$. By  the discussion on the proof of Claim
2.4, the $H$-action on $E_{\Lambda_{\alpha};\g}$ will have empty
fixed point set provided that $H$ can not be conjugated into
$\Lambda_{\alpha}$. From the separability condition (1) and
conjugacy closure condition (2) in the definition of an adapted
family (see Definition 2.1), $H$ can be conjugated into at most one
of the $\Lambda_{\alpha}$; let $k^{-1}$ be the conjugating element,
i.e.\ $k^{-1} H k \subseteq \Lambda_{\alpha}$. Then $H$ fixes
precisely one component $E_{\Lambda_{\alpha};\g}$, namely the
component corresponding to $E_{k \Lambda_{\alpha}}$. Furthermore,
the $H$-action on $E_{k \Lambda_{\alpha}} \cong
E_{\Lambda_{\alpha}}$ is via the $k^{-1} h k$-action on
$E_{\Lambda_{\alpha}}$ which is a model for
$E_{\widetilde{\f}_{\Lambda_{\alpha}}}(\Lambda_{\alpha})$. Since $H
\notin \widetilde{\f}$ , then $k^{-1} H k \notin \widetilde{\f}$,
and this implies $k^{-1} H k \notin
\widetilde{\f}_{\Lambda_{\alpha}}$, therefore $E^H_{k
\Lambda_{\alpha}}=\emptyset$. This immediately gives
$E^H_{\Lambda_{\alpha};\g}=\emptyset$, and hence
$\widehat{X}^H=\emptyset$. This completes the verification of the
second point, and hence the proof of Theorem 2.3.
\end{proof}

\begin{remark}
As we mentioned in Remark 2.2, the collection $\{\g\}$ consisting of
just $\g$ itself is adapted to every pair $(\f, \widetilde{\f})$ of
families of subgroups. Looking at our construction of $\widehat{X}$,
and applying it to $\{\g\}$, we get $\widehat X$ is the join of
$E_{\f}(\g)$ and $E_{\widetilde{\f}}(\g)$. While this is indeed an
$E_{\widetilde{\f}}(\g)$, the construction is not useful, since it
needs an $E_{\widetilde{\f}}(\g)$ to produce an
$E_{\widetilde{\f}}(\g)$. In order to be useful, the construction
requires a ``non-trivial'' adapted family of subgroups.
\end{remark}

\subsection{Relatively hyperbolic groups.}

Our next goal is to exhibit an adapted family of subgroups in the
special case where $\f=\fin$, $\widetilde{\f}=\vc$, and $\g$ is a relatively
hyperbolic group (in the sense of Bowditch \cite{Bo}).  We refer the reader 
to Bowditch for the definition, and will content ourselves with mentioning that
the following classes of groups are relatively hyperbolic: 
\begin{enumerate}
\item free products of finitely many
groups (relative to the factors), 
\item geometrically finite isometry groups of Hadamard
manifolds of pinched negative curvature (relative to maximal parabolic subgroups), 
\item CAT(0)-groups with isolated flats (relative to the flat stabilizers), by recent work of 
Hruska-Kleiner \cite{HK05}, 
\item fundamental groups of spaces obtained via strict relative hyperbolizations (relative 
to the fundamental groups of the subspaces the hyperbolization is taken relative to),
by a recent paper of Belegradek \cite{Bel}.
\end{enumerate}

\begin{theorem}
Let $\g$ be a group which is hyperbolic relative to subgroups
$\{H_i\}^k_{i=1}$, in the sense of Bowditch \cite{Bo} (or
equivalently, a relatively hyperbolic group with the bounded coset penetration property in
the sense of Farb \cite{Fa98}). Consider the collection of subgroups
of $\g$ consisting of:

\begin{enumerate}
\item All conjugates of $H_i$ (these will be called {\it peripheral subgroups}).
\item All maximal virtually infinite cyclic subgroups $V$ such that
$V \nsubseteq g H_i g^{-1}$, for all $i=1, \cdots k$, and for all $g \in \g$.
\end{enumerate}

Then this is an adapted family for the pair $(\fin, \vc)$.
\end{theorem}

\begin{proof}
We first observe that our collection of subgroups clearly satisfies
condition (4) for an adapted family (i.e. every virtually cyclic
subgroup is contained in one of our subgroups).  Furthermore, since
the collection is conjugacy closed by construction, we see that this
collection of subgroups satisfies condition (2) for an adapted
family.  We are left with establishing conditions (1) and (3).

Let us focus on condition (1): any two subgroups in our collection
have finite intersection. A consequence of relative hyperbolicity is
that any two peripheral
subgroups intersect in a finite group (see Section 4 in Bowditch \cite{Bo}). 
Hence to establish (1), it is
sufficient to show:
\begin{itemize}
\item if $V_1,V_2$ are a pair of distinct maximal infinite virtually cyclic
subgroups which {\it do not} lie inside peripheral subgroups, then
$V_1\cap V_2$ is finite, and
\item if $V$ is a maximal infinite virtually cyclic subgroup which
{\it does not} lie in a peripheral subgroup, then it intersects each
$H_i$ in a finite subgroup.
\end{itemize}
To see the first of these two cases, assume that $V_1,V_2$ are as
above, and that $|V_1\cap V_2|=\infty$.  Choose an element $g\in
V_1\cap V_2$ of infinite order, and note that $g$ is {\it
hyperbolic}, in the sense that it {\it cannot} be conjugated into
one of the cusps.  Osin has established (Theorem 4.3 in \cite{Os1})
that the $\mathbb Z$-subgroup $\langle g\rangle$ generated by such
a $g$ lies in a {\it
unique} maximal infinite virtually cyclic subgroup.  This
immediately forces $V_1=V_2$.

Now consider the second of these two cases, and assume $V\cap H_i$
is infinite.  Then picking an element $g \in V\cap H_i$ of infinite
order, we observe that $g$ is of parabolic type (since it lies in
$H_i$).  On the other hand, $V$ contains a finite index subgroup
isomorphic to $\mathbb Z$ generated by an element $h$ of hyperbolic
type.  Now consider the intersection $\langle h\rangle \cap \langle
g\rangle$, and observe that this intersection is non-empty (since
both $\langle h\rangle$ and $\langle g\rangle$ are finite index
subgroups in $V$), hence contains an element which is simultaneously
of hyperbolic type and of parabolic type.  But this is impossible,
giving us a contradiction.

This leaves us with establishing property (3) of an adapted family.
We first note that in a relatively hyperbolic group, the peripheral
subgroups are self-normalizing (see Section 4 in Bowditch 
\cite{Bo}). So we merely need to establish that the maximal infinite
virtually cyclic subgroups $V$ of hyperbolic type are
self-normalizing.  But Osin has established (Theorem 1.5 and
Corollary 1.7 in \cite{Os2}) that if $g\in \Gamma$ is hyperbolic and
has infinite order, and if $V$ is the unique maximal virtually
infinite cyclic subgroup containing $g$, then $V$ has the property
that $V\cap gVg^{-1}$ is {\it finite} for every $g\in \Gamma -V$.
Since $V$ is infinite, this immediately implies that $V=N_\Gamma
(V)$.  This establishes property (3), and completes the proof of the
Theorem.
\end{proof}

%Next, recall that a $\g$-CW-complex is said to be {\it finite} if it has
%finitely many orbits of cells, and that it is said to be of {\it
%finite type} if it has finitely many orbits of cells in each
%dimension.  From Theorem 2.6, we can immediately deduce the
%following finiteness result:

%\begin{corollary}
%Let $\g$ be a group which is hyperbolic relative to subgroups
%$\{H_i\}^k_{i=1}$, and assume that all of the subgroups $H_i$ have
%one of the following properties:
%\begin{enumerate}
%\item finite dimensional classifying spaces $E_{\vc}H_i$.
%\item finite classifying spaces $E_{\vc}H_i$.
%\item classifying spaces $E_{\vc}H_i$ of finite type.
%\end{enumerate}
%Then there is a model for $E_{\vc}\g$ with the same property.
%\end{corollary}

\begin{remark}
(1) Osin has defined a notion of relative hyperbolicity in terms of
relative Dehn functions.  The results we cite in the proof of
Theorem 2.6 make use of his definition.  However, in a previous
paper, Osin has established that for finitely generated groups, his
notion of relative hyperbolicity coincides with Bowditch's
definition (Theorem 1.5 in \cite{Os1}).  As such, his results apply
to the setting in which we are interested.

\noindent (2) We can always view a {\it hyperbolic} group $\g$ as a
group which is hyperbolic relative to the trivial subgroup $\{1\}$.
In this specific case, our construction for the classifying space
$E_{\vc}(\g)$ coincides with the classifying space constructed by
Juan-Pineda and Leary in \cite{JL}.

\noindent (3) In Theorem 2.3, note that it is important that we are using 
Bowditch's notion of relative hyperbolic group.  If we were instead using 
Farb's notion {\it without} the bounded coset penetration property (known
as {\it weak relative hyperbolicity}), the peripheral subgroups are no longer 
self-normalizing.  Indeed, weak relative hyperbolicity is preserved if one replaces
the peripheral subgroups by subgroups of finite index (while in contrast, relative
hyperbolicity is {\it not}).  

\end{remark}

\vskip 10pt

\section[$\gt)$]{The Maximal infinite virtually cyclic subgroups of $\gt$}

Let $\gt$ be the subgroup of $O^+(3,1)$ that preserves the
standard integer lattice $\mathbb Z^{4} \subset \mathbb R^{3,1}$,
that is, $\gt = O^+(3,1) \cap GL(4, \mathbb Z)$.

Since $\gt$ is a subgroup of the discrete group $GL(4, \mathbb
Z)$, it is also a discrete group of $O^+(3,1)$. The
group $\gt$ is a hyperbolic Coxeter, noncocompact, $3$-simplex
reflection group with fundamental domain its defining Coxeter $3$-simplex
$\Delta^3$ (see \cite[pg. 301]{R94}).

The group $\gt$ is part of a nice family of discrete subgroups
of isometries of hyperbolic $n$-space for which the Farrell-Jones
Isomorphism Conjecture in lower algebraic $K$-theory holds, that is
$\ H^{\gt}_n(E_{\vc}(\gt);\mathbb K\mathbb Z^{-\infty}) \cong
K_n(\mathbb Z\gt)$ for $n<2$ (see \cite[Theorem 2.1]{Or04}).  One of our
intentions in this paper is to use this result and our model for
$E_{\vc}(\g)$ constructed in Section 2 to explicitly compute the lower
algebraic $K$-theory of the integral group ring $\mathbb Z\gt$.
In order to accomplish this task,  we must first classify up to
isomorphism the family $\vc$ of all virtually cyclic subgroups of
$\gt$.  At this point we will like to direct the reader to
\cite{Or04} for more information on the relatively hyperbolic groups
$\Gamma_n = O^+(n,1) \cap GL(n+1, \mathbb Z)$, for $n=3,\dots,9$.

We now proceed to classify up to conjugacy all maximal virtually
infinite cyclic subgroups of $\gt$ of hyperbolic type.  The infinite
virtually cyclic subgroups of parabolic type (or cusp groups) are
virtually infinite cyclic subgroups of the cusp group $P4m$,
a 2-dimensional crystallographic group. These groups have already
been classified by Pearson in \cite[Lemma 2.3]{Pe98}.  For the maximal
virtually infinite cyclic subgroups of hyperbolic type, our approach to the
classification problem is {\it geometric}, as opposed to previous approaches
which were algebraic in nature.

\begin{lemma}
Let $Q \leq \gt$ be a infinite virtually cyclic subgroup of $\gt$ of
hyperbolic type.  Then there exist a geodesic $\gamma\subset \mathbb
H^3$ such that $Q \leq \sta_{\gt}(\gamma)$.
\end{lemma}

\begin{proof}
$Q$ is infinite virtually cyclic, hence contains an infinite cyclic
subgroup $H$ of finite index. Since $Q$ is of hyperbolic type, $H$
stabilizes some geodesic $\gamma$ in $\mathbb H^3$. We want to show
that $Q \cdot \, \gamma = \gamma$. Let $g \in Q$ satisfy
$g\cdot\,\gamma \neq \gamma$, and let $g\cdot\,\gamma = \gamma'$.
Note that $gHg^{-1} \leq Q$ stabilizes $\gamma'$. Since $H$ and
$gHg^{-1}$ are both of finite index in $Q$, their intersection is an
infinite cyclic subgroup, call it $K$. Since $\gamma$ and $\gamma'$
are distinct geodesics and $K$ acts by isometries, we get $|K| <
\infty$, contradicting $K \cong \mathbb Z$.
\end{proof}

Note that a subgroup of the type $\sta_{\gt}(\gamma)$ is always
virtually cyclic.  Lemma 3.1 now reduces the problem of classifying
maximal virtually infinity cyclic subgroups of $\gt$ of hyperbolic
type to the more geometric question of finding stabilizers of
geodesics $\gamma \subset \mathbb H^3$.

In order to do this, first we observe that the $\gt$-action on $\mathbb
H^3$ induces a tessellation of $\mathbb H^3$ by copies of the
fundamental domain $\Delta^3$ of $\gt$. This tessellation is
determined by a collection of totally geodesic copies of $\mathbb
H^2$ lying in $\mathbb H^3$ (each of them corresponding to the faces
of the $\Delta^3$) intersecting in a family of geodesics
(corresponding to the 1-skeleton of $\Delta^3$).

Now the stabilizers of the geodesics will depend on the behavior of
the geodesic; more precisely, will depend on the intersection of the
geodesic with the tessellation of $\mathbb H^3$ by copies of the
fundamental domain $\Delta ^3$. Denote by $p:\mathbb H^3 \rightarrow
\mathbb H^3/\gt \cong \Delta^3$ the canonical projection from
$\mathbb H^3$ to the fundamental domain $\Delta^3$.  We first
establish two easy lemmas.

\begin{lemma}
If $\sta_{\gt}(\gamma)$ is infinite, then $p(\gamma) \subset
\Delta^3$ is periodic.
\end{lemma}

\begin{proof}
This follows from the fact that if $\sta_{\gt}(\gamma)$ is virtually
infinite cyclic, then it contains an element of infinite order,
which must act on $\gamma$ by translations. If $g \in \gt$ is this
element, them $p(x)=p(g\cdot \,x) \in \Delta^3$, forcing periodicity
of $p(\gamma)$.
\end{proof}

\begin{lemma}
Let $\gamma \subset \mathbb H^3$ be an arbitrary geodesic, $x \in
\gamma$ an arbitrary point, and $g \in \sta_{\gt}(\gamma) \subset \gt$ an
arbitrary element.  Then we have $p(x)=p(g \cdot \,x)$.
\end{lemma}

Note that Lemma 3.3 is immediate, since any two points (for instance
$x$, and $g \cdot\,x$) have the same image in $\Delta^3$ provided
they differ by an element in $\gt$ (by the definition of the
fundamental domain).

\begin{definition}
For $\gamma$ any geodesic in $\mathbb H^3$, we say that:

\begin{enumerate}
\item $\gamma$ is of {\it type I}, if $\gamma$ is equal to the intersection
of two of the totally geodesic copies of $\mathbb H^2$ inside $\mathbb
H^3$.

\vspace{.2cm}
\item $\gamma$ is of {\it type II}, if $\gamma$ lies entirely within one of
the totally geodesic copies of $\mathbb H^2 \subset \mathbb H^3$, but
not lying in the intersection of two of the $\mathbb H^2$'s.

\vspace{.2cm}
\item $\gamma$ is of {\it type III}, if $\gamma$ does not lie
within one of the totally geodesic $\mathbb H^2 \subset \mathbb H^3$
arising from the tessellation.
\end{enumerate}
\end{definition}

Note that the type of a geodesic can easily be seen in terms of its
image under the projection map $p$. Indeed, geodesics of {\it  type
I} are those for which $p(\gamma)$ lies in the 1-skeleton of
$\Delta^3$, those of {\it type II} have $p(\gamma)$ lying in
$\partial \Delta^3$, but not in the 1-skeleton of $\Delta^3$, and
those of {\it type III} have non-trivial intersection with
$\inti(\Delta^3)$.

We now make the following easy observation: given any geodesic
$\gamma \subset \mathbb H^3$, invariant under the isometric action of a 
Coxeter group $\Gamma$ on $\mathbb H^3$, there exists a short exact sequence:

\[
0 \rightarrow \fix_{\Gamma}(\gamma) \rightarrow \sta_{\Gamma}(\gamma)
\rightarrow \iso_{\Gamma, \gamma}(\mathbb R) \rightarrow 0,
\]

\noindent where $\fix_{\Gamma}(\gamma) \leq  \Gamma$ is the subgroup of
$\Gamma$ that fixes $\gamma$ pointwise, and $\iso_{\Gamma, \gamma}(\mathbb
R)$ is the induced action of $\sta_{\Gamma}(\gamma)$ on $\gamma$
(identified with an isometric copy of $\mathbb R$).  Furthermore we
have:
\begin{enumerate}
\item the group $\iso_{\Gamma, \gamma}(\mathbb R)$, being a discrete 
cocompact subgroup of the isometry group of $\mathbb R$, has to be
isomorphic to $\mathbb Z$ or $D_\infty$.
\item the group $\fix_{\Gamma}(\gamma) \leq  \Gamma $ acts trivially on the geodesic
$\gamma$, hence can be identified with a finite Coxeter group acting on
the unit normal bundle to a point $p\in \gamma$.
\end{enumerate}
In particular, since $\gt$ is a Coxeter group, we can use this short exact sequence
to get an easy description of stabilizers of type II and type III geodesics.

\begin{proposition}
Let $\gamma$ be a geodesic of type III, with $\sta_{\gt}(\gamma)$ infinite
virtually cyclic, them $\sta_{\gt}(\gamma)$ is isomorphic to either
$\mathbb Z$ or $D_{\infty}$.
\end{proposition}

\begin{proof}
Since $p(\gamma) \cap \inti(F) \neq \emptyset$, we
have that $\gamma$ enters the interior of a fundamental domain in
$\mathbb H^3$. If $g \in \gt$ is arbitrary, then the fundamental
domain and its $g$-translate have disjoint interiors; this forces
$\fix_{\gt}(\gamma)=0$.  From the short exact sequence mentioned
above, we immediately obtain $\sta_{\gt}(\gamma) \cong \mathbb
Z$ or $D_{\infty}$, as desired.
\end{proof}

\begin{proposition}
If $\gamma$ is of type II, and $\sta_{\gt}(\gamma)$ is virtually
infinite cyclic, then $\sta_{\gt}(\gamma)$ is either $\mathbb Z
\times \mathbb Z/2$  or $D_{\infty} \times \mathbb Z/2$.
\end{proposition}

\begin{proof}
Let us consider the short exact sequence:

\[
0\rightarrow \fix_{\gt}(\gamma) \rightarrow \sta_{\gt}(\gamma)
\rightarrow \iso_{\gt, \gamma}(\mathbb R) \rightarrow 0
\]

\noindent where $\iso_{\gt,
\gamma}(\mathbb R)$ is either $\mathbb Z$ or $D_{\infty}$. Now let
us focus on $\fix_{\gt}(\gamma)$. Note that since $\gamma$ is of
type II, it lies in one of the totally geodesic $\mathbb H^2 \subset
\mathbb H^3$, but does not  lie in  the intersection of two such
$\mathbb H^2$. We now have that $\fix_{\gt}(\gamma) \cong \mathbb
Z/2$, given by the reflection in the $\mathbb H^2$ containing
$\gamma$. Thus $\sta_{\gt}(\gamma)$ fits into one of the short exact
sequences:

\begin{enumerate}

\item $0 \rightarrow \mathbb Z/2 \rightarrow \sta_{\gt}(\gamma)
\rightarrow \mathbb Z\rightarrow 0$

\item $0 \rightarrow \mathbb Z/2 \rightarrow \sta_{\gt}(\gamma)
\rightarrow D_{\infty} \rightarrow 0$

\end{enumerate}

We next proceed to show that $\sta_{\gt}(\gamma)$ is either
isomorphic to $\mathbb Z \times \mathbb Z/2$ or $D_{\infty} \times
\mathbb Z/2$ according to whether we are in case (1) or (2). Let us
consider case (1), and pick $g \in \sta_{\gt}(\gamma)$ that acts via
translation along $\gamma$. Consider the $\mathbb H^2 \subset
\mathbb H^3$ containing $\gamma$, and note that $g(\mathbb H^2)$ is
another one of the totally geodesic $\mathbb H^2$ containing
$\gamma$. But since $\gamma$ is of type II there is a unique such
totally geodesic $\mathbb H^2$, hence $g(\mathbb H^2)=\mathbb H^2$,
i.e. $g$ leaves the totally geodesic $\mathbb H^2$ invariant. This
immediately implies that the subgroup generated by $g$ commutes with
the reflection in the $\mathbb H^2$, yielding that
$\sta_{\gt}(\gamma) \cong \mathbb Z \times \mathbb Z/2$, as desired.
The argument in case (2) is nearly identical; pick a pair $g, h \in
\sta_{\gt}(\gamma)$ whose image are the generators for $D_{\infty}$,
and observe that $g, h$ must be reflections in a pair of $\mathbb
H^2$'s both of which are perpendicular to $\gamma$. This implies
that the hyperplanes are both orthogonal to the totally geodesic
$\mathbb H^2$ containing $\gamma$, and hence that $g, h$ both map
the $\mathbb H^2$ to itself. This yields that $g, h$ both commute
with reflection in the corresponding $\mathbb H^2$, and hence that
$\sta_{\gt}(\gamma) \cong D_{\infty} \times \mathbb Z/2$.
\end{proof}

To finish the classification of maximal virtually  infinite cyclic
subgroup of $\gt$, we are left with the study of stabilizers of
geodesic of type I.  Recall that geodesics  $\gamma$ of type I are
those  for which $p(\gamma)$ lies in the 1-skeleton of $\Delta^3$.  In this
situation, the short exact sequence is not too useful for our purposes.  Instead,
we switch our viewpoint, and appeal instead to Bass-Serre theory \cite{S80}.

\begin{figure}[htbp]
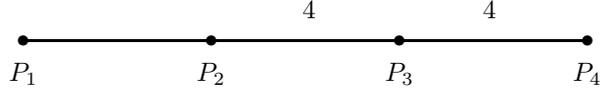

%\begin{center}
\hfil\vbox{\beginpicture
\setcoordinatesystem units <1cm,1cm> point at -1 .8
\setplotarea x from -1 to 9.5, y from -.8 to .8
\linethickness=.8pt
\putrule from 0 0  to 7.5 0 
\put {$\bullet$} at  0  0
\put {$\bullet$} at  2.5 0
\put {$\bullet$} at  5  0
\put {$\bullet$} at 7.5  0
\put {$P_1$} [t] at   0 -.3 
\put {$P_2$} [t] at   2.5 -.3 
\put {$P_3$} [t] at   5 -.3 
\put {$P_4$} [t] at   7.5 -.3 
\put {4} [b] at   3.8 .3 
\put {4} [b] at   6.2 .3 
\endpicture}\hfil
%\psfrag{s1}[c]{$P_1$}
%\psfrag{s2}[c]{$P_2$}
%\psfrag{s3}[c]{$P_3$}
%\psfrag{s4}[c]{$P_4$}
%\psfrag{4}[c]{$4$}
%\psfrag{4}[c]{$4$}
%\includegraphics{coxeter1.eps}
%\end{center}
\caption{The Coxeter graph of $\Gamma_3=[3,4,4]$}
\end{figure}

Recall that the group $\Gamma_3$ is a Coxeter group with Coxeter graph given in Figure 1.
The fundamental domain for the $\gt$-action on $\mathbb H^3$ is a 3-simplex $\Delta ^3$ in $\mathbb H^3$,
with one ideal vertex.  The group $\gt \leq Isom(\mathbb H^3)$ is generated by a reflection in the
hyperplanes (totally geodesic $\mathbb H^3$) extending the four faces of the 3-simplex.  Each 
reflection corresponds to a generator, and the faces of the 3-simplex are labelled $P_1, \ldots ,P_4$
according to the corresponding generator.  Note that the angles between the faces can be read off
from the Coxeter graph.  The possible type I geodesics correspond to geodesics extending 
the intersections of pairs of faces (and hence, there are at most six such geodesics).

\begin{proposition}
Let $\gamma$ be the geodesic extending  the intersection of the two
hyperplanes $P_1 \cap P_3$, and let $\eta$ be the geodesic extending
$P_1 \cap P_2$. Then $\sta_{\gt}(\gamma) \cong D_2 \times
D_{\infty}$, and $\sta_{\gt}(\eta) \cong D_3 \times D_{\infty}$.
Furthermore, for the geodesics extending the remaining edges in the
one skeleton, the stabilizers are finite.
\end{proposition}

\begin{proof}
Our procedure for identifying the stabilizers of the geodesics of
type I relies on the observation that when these stabilizers act
cocompactly on the geodesics, Bass-Serre theory gives us an easy
description of the corresponding stabilizer. Indeed, if the quotient
of the geodesic in the fundamental domain is a segment, then the
stabilizer of the geodesic is an amalgamation of the vertex
stabilizers, amalgamated over the edge stabilizer. Let us carry out
this procedure in the specific case of $\Gamma _3$.  Notice that
three of the six edges in the 1-skeleton of the fundamental domain
$\Delta ^3$ are actually geodesic rays (since one of the vertices is
an ideal vertex), and hence their geodesic extensions will have
finite stabilizers.  This leaves us with three edges in the
1-skeleton to worry about, namely those corresponding to the edges
$P_1\cap P_2$, $P_1\cap P_3$, and $P_1\cap P_4$.  Let us first focus on 
the edge
corresponding to the intersection $P_1\cap P_4$.  We claim that the
geodesic extending this edge projects into a non-compact segment
inside $\Delta ^3$.  Indeed, if one considers the link of the vertex
$P_1\cap P_2\cap P_4$ (a 2-dimensional sphere), we note that it has
a natural action by a parabolic Coxeter group, namely the stabilizer
of this vertex.  Corresponding to this action is a tessellation of
$S^2$ by geodesic triangles, where the vertices of the triangles
correspond to the directed edges from the given vertex to
(translates) of the remaining vertices from the fundamental domain.
In other words, each vertex of the tessellation of $S^2$ comes
equipped with a {\it label} identifying which directed edge in the
fundamental domain it corresponds to.  But now in the tessellation of
the vertex $P_1\cap P_2\cap P_4$, it is easy to see that the vertex
antipodal to the one corresponding to the edge $P_1\cap P_4$, has a
label indicating that it projects to the edge $P_2\cap P_4$.  This
immediately allows us to see that the geodesic ray extending the
edge $P_1\cap P_4$ does {\it not} project to a closed loop in
$\Delta ^3$, and hence has finite stabilizer.

%% Test for figure 1 10/12/05
%%%%%%%%%%%%%%%%%%%%%%%%%%% Figure 1%%%%%%%%%%%%%%%%%%%%%%%%%%%%%%%%%%
\begin{figure}[htbp]
\begin{center}
\hfil\vbox{\beginpicture
\setcoordinatesystem units <.7cm,.7cm> point at -1 1
\setplotarea x from -.2 to 14, y from -2 to 1
\setlinear
\linethickness=.8pt
%\linethickness=1.2pt  %Was used for ``big dots''
\putrule from 0 0 to 4 0
\putrule from 8 0 to 12 0
%\put {\circle*{2}} at 5.8 0
\put {$\bullet$} at  0  0
\put {$\bullet$} at  4 0
\put {$\bullet$} at  8 0
\put {$\bullet$} at  12 0
\put {$\overline{\eta}$}  at 2 1.3
\put {$\overline{\gamma}$} at 10 1.3
%put {or} at 6 0
\put {$D_6$} [t] at   0 -.2
\put {$D_3$} [t] at   2 -.2
\put {$D_6$} [t] at   4 -.2
\put {$D_2\times Z/2$} [t] at   8 -.2
\put {$D_2$} [t] at   10 -.2
\put {$D_2\times Z/2$} [t] at   12 -.2
%\put {Figure 1} [t] at 6 -1.5
\endpicture}\hfil
\end{center}
\caption{}
\end{figure}
%%%%%%%%%%%%%%%%%%%%%%%%%%%%%%%%%%%%%%%%%%%%%%%%%%%%%%%%%%%%%%%%%%%%%%%

So we are left with identifying the stabilizers of the two geodesics
$\gamma$ and $\eta$. Let us consider the geodesic $\gamma$, and
observe that the segment that is being extended joins the vertex
$P_1\cap P_2\cap P_3$ to the vertex $P_1\cap P_3\cap P_4$. When
looking at the tilings of $S^2$ one obtains at each of these two
vertices, we find that in both cases, the label of the vertex
corresponding to the edge $P_1\cap P_3$ is antipodal to a vertex
with the same labeling.  This tells us that the geodesic $\gamma$
projects precisely to the segment $P_1\cap P_3$ in the fundamental
domain $\Delta ^3$.  Now the subgroup of $Stab_{\gt}(\gamma)$ that
fixes the vertex $P_1\cap P_2\cap P_3$ is the subgroup of the
Coxeter group stabilizing the vertex, that additionally fixes the
pair of antipodal vertices in the tiling of $S^2$ corresponding to
$\gamma$.  But it is immediate by looking at the tiling that this
subgroup at each of the two vertices is just $D_2\times \mathbb
Z_2$, where in both cases the $D_2$ stabilizer of the edge $P_1\cap
P_3$ injects into the first factor, and the $\mathbb Z_2$-factor is
generated by the reflection of $S^2$ that interchanges the given
antipodal points.  This yields that $Stab_{\gt} (\gamma)$ is the
amalgamation of two copies of $D_2\times \mathbb Z_2$ along the
$D_2$ factors, giving us $D_2\times D_\infty$ (see Figure 2 for the
graph of groups). An identical analysis
in the case of the geodesic $\eta$ yields that the stabilizer is an
amalgamation of two copies of $D_3\times \mathbb Z_2$ along the
$D_3$ factors, yielding that the stabilizer is $D_3\times D_\infty$ (see Figure 2).   
This completes the proof of Proposition 3.8.
\end{proof}

\section{The computations of  $Wh_n(\gt), n<2$}

In this section we briefly recall the Isomorphism Conjecture in lower algebraic $K$-theory (the interested reader should refer to \cite{FJ93}, \cite{DL98}), and compute the homology groups $H^{\gt}_n(E_{\vc}(\gt); \mathbb K\mathbb Z^{-\infty}) \cong Wh_n(\gt)$ for $n<2$.

The {\it Farrell and Jones Isomorphism Conjecture in algebraic
$K$-theory}, reformulated in terms of the Davis and L\"uck functor $\mathbb KR^{-\infty}$ (see \cite{DL98}), states that the assembly map 
 $H^{\Gamma}_n(E_{\vc}(\g);\mathbb KR^{-\infty})  
\longrightarrow K_n(R\Gamma)$
is an isomorphism for all $n \in \mathbb Z$.

The main point of the validity of this conjecture is that it allows
the computations of the groups of interest
$K_n(R\Gamma)$ from the  values of $\mathbb KR^{-\infty}(\Gamma/H)$ on
the groups  $H \in \mathcal V\mathcal C$.

The pseudo-isotopic  version of the Farrell and
Jones Conjecture is obtained by replacing the algebraic $K$-theory spectrum
by the functors $\mathcal P_\ast$, $\mathcal
P^{\text{diff}}_\ast$, which map from the category of topological spaces X
to the category of $\Omega-\mathcal S\mathcal P\mathcal E\mathcal
C\mathcal T\mathcal R\mathcal A$. The functor $\mathcal P_\ast(?)$
(or $\mathcal P^{\text{diff}}_\ast(?))$ maps the space $X$ to the
$\Omega$-spectrum of stable topological (or smooth) pseudo-isotopies
of $X$ (see \cite[Section 1.1]{FJ93}). 

The relation between $\mathcal P_\ast(? )$ and lower algebraic $K$-theory
is given by the work of Anderson and Hsiang \cite[Theorem 3]{AH}. They show 
\[
\pi_j(\mathcal P_\ast(X))=
\begin{cases}
Wh(\mathbb Z\pi_1(X)), & j=-1 \\
\tilde {K}_0(\mathbb Z\pi_1(X) ), & j=-2 \\
K_{j+2}(\mathbb Z \pi_1(X)), & j \leq -3.
\end{cases}
\]
The main result in \cite{FJ93} is that the Isomorphism Conjecture is
true for the pseudo-isotopy and smooth pseudo-isotopy functors when
$\pi_1(X)=\Gamma$ is a discrete cocompact subgroup of a
virtually connected Lie group. This result together with
the identification given by Anderson and Hsiang of the lower homotopy
groups of the pseudo-isotopy spectrum and the lower algebraic
$K$-theory implies the following  Theorem (see \cite[Section 1.6.5, and
Theorem 2.1]{FJ93}): 

\begin{theorem}[Farrell, F.T. and Jones, L.E]
Let $\Gamma$ be a cocompact discrete subgroup of a virtually connected Lie
 group. Then the assembly map
\[
 H_n^{\g}(E_{\vc}(\g);\mathbb K\mathbb Z^{-\infty})\longrightarrow  K_n(\mathbb Z\g) 
\]
is an isomorphism for $n \leq 1$ and a surjection for $n=2$.
\end{theorem}

Farrell and Jones also proved Theorem 4.1 for discrete cocompact
groups, acting properly discontinuously by isometries on a simply
connected Riemannian manifold $M$ with everywhere non-positive
curvature (\cite[Proposition 2.3]{FJ93}). Berkove, Farrell, Pineda, and Pearson
extended this result to  discrete groups, acting properly discontinuously on
hyperbolic  $n$-spaces via isometries, whose orbit space has finite
volume (but not necessarily compact), (see \cite[Theorem A]{BFPP00}). In
particular this result is valid for $\Gamma$  a hyperbolic, non-cocompact,
$n$-simplex reflection group. 

\vspace{.2cm}
Therefore for  $\g=\gt$, it follows that
\[
H_n^{\gt}(E_{\vc}(\gt);\mathbb K\mathbb Z^{-\infty}) \cong Wh_n(\mathbb Z\gt), \quad \text{ for} \quad n<2. 
\]

\vspace{.2cm}
Hence to compute the lower algebraic $K$-theory of the integral group ring $\gt$ it suffices to compute for $n<2$, the homotopy groups

\[
H_n^{\gt}(E_{\vc}(\gt);\mathbb K\mathbb Z^{-\infty}).
\]

\vspace{.2cm}
These computations are feasible using the Atiyah-Hirzebruch type spectral
sequence described by  Quinn in  \cite[Theorem 8.7]{Qu82} for the pseudo-isotopy spectrum $\mathcal P  $ (see Section 5): 

\[
E^2_{p,q}=H_p(E_{\f}(\g)/\g\; ;\{Wh_q(\mathbb Z \g_{\sigma})\})
\Longrightarrow 
Wh_{p+q}(\g),
\]

\noindent
where
\[
Wh_q(F)=
\begin{cases}
Wh(F), & q=1 \\
\tilde {K}_0(\mathbb Z F), & q=0 \\
K_n(\mathbb Z F), & q \leq -1.
\end{cases}
\]

All the information need it to compute the $E^2$ term is encoded in $E_{\vc}(\gt)/\gt$ and the algebraic $K$-groups of the finite subgroups and maximal infinite virtually cyclic subgroups of $\gt$.

%%%%%%%%%%%%%%%%%%%%%%%%%%%%%%%%%%%%%%%%%%%%%%%%%%%%%%%%%%

We now proceed to give a proof of our main Theorem 1.1.  Recall that this Theorem 
states that the lower algebraic $K$-theory of the integral group ring of $\gt$ is
given as follows:
\[
\begin{aligned}
Wh(\gt) &\cong \bigoplus_{\infty} \mathbb Z/2\\
\tilde{K}_0(\mathbb Z \gt) &\cong
\mathbb Z/4 \oplus \mathbb Z/4 \oplus \bigoplus_{\infty} \mathbb Z/2\\
K_{-1}(\mathbb Z \gt) &\cong \mathbb Z \oplus \mathbb Z, \quad \text{and}\\
K_n(\mathbb Z \gt) &\cong 0, \quad \text{for} \; n < -1.
\end{aligned}
\]

\begin{proof}

Since $\gt$ is hyperbolic relative to $P4m$, then the
fundamental domain $\Delta^3$ has one cusp with cusp subgroup $P4m$
(a 2-dimensional crystallographic group), which is the unique
maximal parabolic subgroup of $\gt$ (up to conjugacy).

Let $X_0 \subset \mathbb H^3$ be the space obtained by truncating the
cusp, that is, we remove from $\mathbb H^3$ a countable collection
of (open) horoballs, in a $\gt$-equivariant way, where the
horoballs $B_i$ are based at the points $\{p_i\} \subset
\partial^{\infty}\mathbb H^3=S^2$ which are fixed points of the
subgroups $G \leq \gt$ which are conjugate to $P4m$.  Now let $X=
X_0/\gt$, i.e.\  $X=\Delta^3 \setminus B_i$, where $B_i$ is the open horoball
from our countable collection that is based at the ideal vertex of $\Delta^3$.

Note that $X$ satisfies the requirements to be a model for
$E_{\fin}(\gt)/\gt$. $X$ has five faces (with stabilizers: 1,
$\mathbb Z/2$), nine edges (with stabilizers: $\mathbb Z/2$, $ D_2$,
$ D_3$, $ D_4$), and six vertices (with stabilizers: $D_2$, $D_4$,
$D_6$, $\mathbb Z/2 \times D_4$, $\mathbb Z/2 \times S_4$).

Let $\vcm$ be the family of all maximal virtually cyclic subgroups
of $\gt$ of hyperbolic type.  Form the space $Z$ consisting of one
point, denoted $\lbrack gV \rbrack$ for each left coset of $V \in \vcm$, i.e.\
$Z=\coprod_{V \in \vcm, g\in \gt} \lbrack gV \rbrack$. Let $E_{\vc}(P4m)$ the
classifying space for $P4m$ constructed by Alves and Ontaneda in
\cite{AO03} ($E_{\vc}(P4m)$ is a  4-dimensional CW-complex, for the
isotropy groups we refer the reader to \cite{AO03}), and let $Y= Z
\coprod E_{\vc}(P4m)$. Form the space:

\[
\widehat{X}=X \ast Y = X \ast \bigg \lbrack \bigg(\coprod_{V \in
\vcm, g\in \gt} \lbrack gV \rbrack \bigg) \coprod \bigg(\coprod_{g\in \gt} E_{\vc}(P4m)  \bigg) \bigg\rbrack
\]

\noindent From Theorem 2.3 and Theorem 2.6, we have that
$\widehat{X}$ satisfies the requirements to be a model for
$E_{\vc}(\gt)/\gt$. $\widehat{X}$ is an 8 dimensional
$\gt$-CW complex with isotropy groups:

\begin{enumerate}

\item {\bf Isotropy of the 0-cells:}
The stabilizers of the 0-cells are either finite or maximal infinite
virtually cyclic of hyperbolic type:

\vspace{.2cm}
\begin{enumerate}
\item {\it Finite subgroups:} $\mathbb Z/2$, $D_2$, $D_4$, $D_6$,
$\mathbb Z/2 \times D_4$, $\mathbb Z/2 \times S_4$.

\vspace{.2cm}
\item {\it Maximal infinite virtually cyclic subgroups of hyperbolic type:}
$\mathbb Z$, $D_{\infty}$, $\mathbb Z \times \mathbb Z/2$,
$D_{\infty} \times \mathbb Z/2$,  $D_2 \times D_{\infty} \cong D_2
\times \mathbb Z/2 \ast_{D_2} D_2 \times \mathbb Z/2$, $D_3 \times
D_{\infty} \cong D_6 \ast_{D_3} D_6$.
\end{enumerate}

\vspace{.3cm}
\item {\bf Isotropy of the 1-cells:}
The stabilizers of the 1-cells are either finite or infinite
virtually cyclic of parabolic type:

\begin{enumerate}
\vspace{.2cm}
\item {\it Finite subgroups:} 1, $\mathbb Z/2$, $D_2$, $D_3$, $D_4$, $D_6$.

\vspace{.2cm}
\item {\it Infinite virtually cyclic subgroups of parabolic type:}
$\mathbb Z$, $D_{\infty}$, $\mathbb Z \times \mathbb Z/2$,
$D_{\infty} \times \mathbb Z/2$.
\end{enumerate}

\vspace{.3cm}
\item {\bf Isotropy of the 2-cells:} The stabilizers of the 2-cells are
either finite or infinite virtually cyclic of parabolic type:

\vspace{.2cm}
\begin{enumerate}
\item {\it Finite subgroups:} 1, $\mathbb Z/2$, $D_2$, $D_3$, $D_4$.

\vspace{.2cm}
\item {\it Infinite virtually cyclic subgroups of parabolic type:} $\mathbb Z$, $D_{\infty}$
\end{enumerate}

\vspace{.3cm}
\item {\bf Isotropy of the 3-cells and 4-cells:} 1, $\mathbb Z/2$.

\vspace{.3cm}
\item {\bf Isotropy of the 5-cells, 6-cells, 7-cells and 8-cells:} Trivial.
\end{enumerate}

\vspace{.2cm} The complex that gives the homology of
$E_{\vc}(\gt) /\gt$ with local coefficients $\{Wh_q(
F_{\sigma})\}$ has the form
\[
\bigoplus_{{\sigma}^8}^{}Wh_q( F_{{\sigma}^8})\rightarrow \cdots \rightarrow
\bigoplus_{{\sigma}^2}^{}Wh_q( F_{{\sigma}^2}) \rightarrow
\bigoplus_{{\sigma}^1}^{}Wh_q( F_{{\sigma}^1}) \rightarrow
\bigoplus_{{\sigma}^0}^{}Wh_q( F_{{\sigma}^0}),
\]
where ${\sigma}^i$ denotes the cells in dimension $i$, and $Wh_q(
F_{{\sigma}^i})$ occurs in the summand as many times as the numbers
of conjugacy classes of the subgroup $F_{{\sigma}^i}$ in $\gt$.
The homology of this complex will give us the entries for the
$E^2$-term of the spectral sequence.  We now proceed to 
analyze this complex for each of the following cases: $q<-1$,
$q=-1,0,1$.

\vspace{.2cm} \noindent {$\mathbf {q<-1}$}. Carter showed in \cite{C80} that
$K_q(\mathbb Z F)=0$ when $F$ is a finite group. Farrell and Jones
showed in \cite{FJ95}  that $K_q(\mathbb Z Q)=0$ when $Q$ is a
infinite virtually cyclic group. Hence the whole complex consists of
zero terms and we obtain  $E^2_{p,q}=0$ for $q<-1$.

\vspace{.2cm} \noindent {$\mathbf {q=-1}$}. Again using Carter's result in
\cite{C80}, $K_{-1}(\mathbb Z F)=0$, for all the finite subgroups
which occur as stabilizers of the $n$-cells, with $n=2, \cdots, 8$.
For $n=2$, we have 2-cells with stabilizers which are virtually
infinite cyclic subgroups of parabolic type. Bass in \cite{Bas68}
showed that $K_n(\mathbb ZQ)$ for $n<0$ vanishes if $Q=\mathbb Z$, or
$Q=D_{\infty}$. Therefore $E^2_{p,-1}=0$ for $p \geq 2$. Also
$K_{-1}(\mathbb Z F)=0$, for all finite subgroups which are
stabilizers of the 1-cells except for $F=D_6$; we have two 1-cells
with stabilizer $F=D_6$, for which $K_{-1}(\mathbb Z D_6)=\mathbb Z$
(see \cite[pg. 274]{Pe98}). For $n=1$, we also have 1-cells with
stabilizers which are infinite virtually cyclic subgroups of
parabolic type. Pearson in \cite{Pe98} showed that
$K_n(\mathbb ZQ)$ for $n<0$ vanishes if $Q=\mathbb Z \times \mathbb
Z/2$ or  $Q=D_{\infty} \times \mathbb Z/2$. It follows that for
$p=0, 1$ the complex  may have non-zero terms in dimension one and
zero:

\[
\cdots \rightarrow 0 \rightarrow \bigoplus_{{\sigma}^1}^{}K_{-1}(
D_6) \rightarrow \bigg\lbrack \thinspace \negthickspace
\bigoplus_{F_{\sigma^0} \in \fin}  \negthickspace \negthickspace
\negthickspace K_{-1}(\mathbb Z F_{{\sigma}^0}) \thinspace
\thinspace \thinspace \oplus \bigoplus_{Q_{\sigma^0} \in
\vc_{\infty}} \negthickspace \negthickspace \negthickspace
K_{-1}(\mathbb Z Q_{{\sigma}^0})\bigg\rbrack
\]

If $F_{\sigma^0}$ is one of the finite subgroups groups $\mathbb
Z/2$, $D_2$, or $D_4$,  then $K_{-1}(\mathbb ZF)=0$ (see
\cite{C80}). As we mentioned earlier,  $K_{-1}(\mathbb ZD_6)= \mathbb
Z$ (see \cite{Pe98}). If   $F= \mathbb Z/2 \times D_4$, or $F=
\mathbb Z/2 \times S_4$,  then Ortiz in \cite[pg. 350]{Or04}) showed
that  $K_{-1}(\mathbb Z[\mathbb Z/2 \times D_4])=0$, and
$K_{-1}(\mathbb Z[\mathbb Z/2 \times S_4])=\mathbb Z$.

If $Q_{\sigma^0} \in \vc_{\infty}$, then $Q_{\sigma^0}$ is one of
the groups: $\mathbb Z$, $D_{\infty}$, $\mathbb Z \times \mathbb
Z/2$, $D_{\infty} \times \mathbb Z/2$,  $D_2 \times D_{\infty} \cong
D_2 \times \mathbb Z/2 \ast_{D_2} D_2 \times \mathbb Z/2$, $D_3
\times D_{\infty} \cong D_6 \ast_{D_3} D_6$. As we mentioned earlier,
the groups $K_n(\mathbb ZQ)$ for $n<0$ vanish if $Q=\mathbb Z$, or
$Q=D_{\infty}$  (see \cite{Bas68}). For the groups with $\mathbb
Z/2$ summands, we also have that $K_{n}(\mathbb ZQ) = 0$ if  $n <0$,
(see \cite[pg. 272]{Pe98}).

The other two cases are groups of the form $Q=Q_0\ast_{D_n}
Q_1$, with $n=2,3$.  In \cite{FJ95} Farrell and Jones show that if $Q$ is infinite
virtually cyclic, then $K_n(\mathbb Z Q)$ is zero for $n <-1$ and
that $K_{-1}(\mathbb Z Q)$ is generated by the images of
$K_{-1}(\mathbb Z F)$ where $F$ ranges over all finite subgroups $F
\subset Q$. Since $Wh_q(D_2)=0$ for $q<0$, and $Wh_q(D_2 \times
\mathbb Z/2)=0$ for $q<0$ (see \cite{LS00}),  then for  $Q=D_2
\times D_{\infty} \cong D_2 \times \mathbb Z/2 \ast_{D_2} D_2 \times
\mathbb Z/2$, it  follows that $K_{-1}(\mathbb ZQ)=0$.  Since
$K_{-1}(\mathbb ZD_3)=0$, and $K_{-1}(\mathbb ZD_6)=\mathbb Z$ (see
\cite{Or04}, \cite{Pe98}), then  for $Q=D_3 \times D_{\infty} \cong
D_6 \ast_{D_3} D_6$, it  follows that $K_{-1}(\mathbb ZQ)=\mathbb Z
\oplus \mathbb Z$.

Since there is only one conjugacy class for each of the subgroups of
$\gt$ occurring as stabilizers of the 1-cells and the 0-cells, then
the complex that gives the homology groups
$H_0(\widehat{X};\{K_{-1}(\mathbb Z F_{\sigma})\})$, and
$H_1(\widehat{X};\{K_{-1}(\mathbb Z F_{\sigma})\})$ yields the
following exact sequence
\[
0 \rightarrow \mathbb Z \oplus \mathbb Z \rightarrow \mathbb Z
\oplus \mathbb Z \oplus \mathbb Z \oplus \mathbb Z \rightarrow 0.
\]
Hence after working through the exact sequence, we have that that
the $E^2$ term for $\widehat{X}$ has the following entries for
$q=-1$: $E^2_{p, q}=0$ for $p \geq 1$, and $E^2_{0,-1}=\mathbb Z
\oplus \mathbb Z$.

\vspace{.2cm} \noindent {$\mathbf {q=0}$}. It is well known that $\tilde
{K}_0(\mathbb Z F)=0$ when $F$ is any of the finite subgroups that
occur as stabilizers of the n-cells, for $n=1, \cdots,8$ (see for
example \cite{Re76}, \cite{Ro94}). For $n=1, 2$,  we have 1-cells
and 2-cells with stabilizers which are infinite virtually cyclic
subgroups of parabolic type. In \cite{Pe98} Pearson showed that
$\tilde{K}_0(\mathbb Z Q)=0$ if $Q=\mathbb Z$, $D_{\infty}$,
$\mathbb Z \times \mathbb Z/2$, $D_{\infty} \times \mathbb Z/2$,
therefore $E^2_{p,0}=0$ for $p \geq 1$. For $p=0$ the complex may
have non-zero terms in dimension zero and the resulting homology is:
\[
H_0(\widehat{X} ;\{\tilde{K}_0(\mathbb Z F_{\sigma})\}) =
\negthickspace \bigoplus_{F_{\sigma^0} \in \fin} \negthickspace
\negthickspace \negthickspace\tilde{K}_{0}(\mathbb Z
F_{{\sigma}^0})\thinspace \thinspace \thinspace \oplus
\bigoplus_{Q_{\sigma^0} \in \vc_{\infty}}  \negthickspace
\negthickspace \negthickspace \tilde{K}_{0}(\mathbb Z
Q_{{\sigma}^0}).
\]

If $F_{\sigma^0}$ is one of the finite subgroups $\mathbb Z/2$,
$D_2$, $D_4$, or $D_6$, then $\tilde{K}_{0}(\mathbb ZF)=0$ (see
\cite{Re76}). If $F= \mathbb Z/2 \times D_4$, or  $F= \mathbb Z/2
\times S_4$, Ortiz in \cite[pg. 351]{Or04} showed that
$\tilde{K}_{0}(\mathbb Z[\mathbb Z/2 \times D_4]) =\mathbb Z/4$, and
$\tilde{K}_{0}(\mathbb Z[\mathbb Z/2 \times S_4])=\mathbb Z/4$.

If $Q_{\sigma^0}$  is one of the maximal infinite virtually cyclic
subgroups $\mathbb Z$, $D_{\infty}$, $\mathbb Z \times \mathbb Z/2$,
or $D_{\infty} \times \mathbb Z/2$,  then $\tilde{K}_0(\mathbb Z
Q)=0$ (see \cite{Pe98}).  For the remaining subgroups, using
\cite[Lemma 3.8]{CP02}, we have that for $Q=D_2 \times \mathbb Z/2
\ast_{D_2}  D_2 \times \mathbb Z/2$, $\tilde {K}_0(\mathbb ZQ) \cong
NK_0(\mathbb ZD_2; B_1, B_2)$, where $B_i=\mathbb Z[D_2 \times Z/2
\setminus D_2]$ is the $\mathbb ZD_2$ bi-module generated by $D_2
\times \mathbb Z/2 \setminus D_2$ for $i=1,2$. For $Q=D_6 \ast_{D_3}
D_6$, we have that $\tilde {K}_0(\mathbb ZQ) \cong NK_0(\mathbb
ZD_3; C_1, C_2)$, where $C_i=\mathbb Z[D_6 \setminus D_3]$ is the
$\mathbb ZD_3$-bimodule generated by $D_6 \setminus D_3$ for
$i=1,2$. The Nil-groups $NK_0$ appearing in these computations are
the Waldhausen's Nil-groups.

Hence we obtain that
\[
E^2_{0,0}= \mathbb Z/4 \oplus \mathbb Z/4 \oplus NK_0(\mathbb ZD_2;B_1, B_2)  \oplus NK_0(\mathbb ZD_3;C_1, C_2).
\]

\vspace{.2cm} \noindent $\mathbf {q=1}$. Oliver in \cite{O89} showed that
$Wh(F)=0$ when $F$ is any of the finite subgroups that occur as
stabilizers of the $n$-cells, for $n=1, \cdots, 8$. For $n=1, 2$,
we have 1-cells and 2-cells with stabilizers which are virtually
infinite cyclic subgroups of parabolic type. In \cite{Pe98} Pearson
show that $Wh(Q)=0$ if $Q=\mathbb Z$, $D_{\infty}$, $\mathbb Z
\times \mathbb Z/2$, $D_{\infty} \times \mathbb Z/2$, therefore
$E^2_{p,1}=0$ for $p \geq 1$. As before for $p=0$ the complex may
have non-zero terms in dimension zero and the resulting homology is:
\[
H_0(\widehat{X} ;\{Wh(F_{\sigma})\}) = \negthickspace
\bigoplus_{F_{\sigma^0} \in \fin} \negthickspace \negthickspace
\negthickspace Wh(F_{{\sigma}^0})\thinspace \thinspace \thinspace
\oplus \bigoplus_{Q_{\sigma^0} \in \vc_{\infty}}  \negthickspace
\negthickspace \negthickspace Wh(Q_{{\sigma}^0}).
\]

If $F_{\sigma^0}$ is one of the finite subgroups $\mathbb Z/2$,
$D_2$, $D_4$, or $D_6$, then $Wh(F)=0$ (see \cite{O89}). If $F=
\mathbb Z/2 \times D_4$, or  $F= \mathbb Z/2 \times S_4$, Ortiz in
\cite[pg. 352]{Or04} showed that $Wh(\mathbb Z/2 \times D_4) =0$, and
$Wh(\mathbb Z/2 \times S_4)=0$.

If $Q_{\sigma^0}$  is one of the maximal infinite virtually cyclic
subgroups $\mathbb Z$, $D_{\infty}$, $\mathbb Z \times \mathbb Z/2$,
or $D_{\infty} \times \mathbb Z/2$,  then $Wh(Q)=0$ (see
\cite{Pe98}).  For the remaining subgroups, using \cite[Lemma
3.8]{CP02}, we have that for $Q=D_2 \times \mathbb Z/2 \ast_{D_2}
D_2 \times \mathbb Z/2$, $Wh(Q) \cong NK_1(\mathbb ZD_2; B_1, B_2)$,
where $B_i=\mathbb Z[D_2 \times Z/2 \setminus D_2]$ is the $\mathbb
ZD_2$ bi-module generated by $D_2 \times \mathbb Z/2 \setminus D_2$
for $i=1,2$. For $Q=D_6 \ast_{D_3} D_6$, we have that $Wh(Q) \cong
NK_1(\mathbb ZD_3; C_1, C_2)$, where $C_i=\mathbb Z[D_6 \setminus
D_3]$ is the $\mathbb ZD_3$-bimodule generated by $D_6 \setminus
D_3$ for $i=1,2$. The Nil-groups $NK_1$ appearing in these
computations are the Waldhausen's Nil-groups. It follows that
\[
E^2_{0,1}= NK_1(\mathbb ZD_2;B_1, B_2)  \oplus NK_1(\mathbb ZD_3;C_1, C_2).
\]

Hence the spectral sequence collapses at $E^2$, giving the following preliminary results on the algebraic $K$-groups $Wh_n(\gt)$ for $n<2$.

\[
\begin{aligned}
Wh(\gt) &\cong NK_1(\mathbb ZD_2;B_1, B_2)  \oplus NK_1(\mathbb ZD_3;C_1, C_2)\\
\tilde{K}_0(\mathbb Z \gt) &\cong
\mathbb Z/4 \oplus \mathbb Z/4 \oplus NK_0(\mathbb ZD_2;B_1, B_2)  \oplus NK_0(\mathbb ZD_3;C_1, C_2)\\
K_{-1}(\mathbb Z \gt) &\cong \mathbb Z \oplus \mathbb Z, \; \text{and}\\
K_n(\mathbb Z \gt) &\cong 0, \; \text{for} \; n < -1.
\end{aligned}
\]
\noindent where $B_i=\mathbb Z[D_2 \times \mathbb Z/2 \setminus D_2]$ is the $\mathbb
ZD_2$ bi-module generated by $D_2 \times \mathbb Z/2 \setminus  D_2$
for $i=1,2$, and $C_i=\mathbb Z[D_6 \setminus D_3]$ is the $\mathbb
ZD_3$-bimodule generated by $D_6 \setminus D_3$ for $i=1,2$. 

Finally, to conclude the proof, we observe that the Waldhausen Nil-groups that appear above
are isomorphic to:
$$NK_1(\mathbb ZD_2;B_1, B_2) \cong NK_0(\mathbb ZD_2;B_1, B_2) \cong \bigoplus_\infty \mathbb Z/2$$
$$NK_1(\mathbb ZD_3;C_1, C_2) \cong NK_0(\mathbb ZD_3;C_1, C_2) \cong 0$$
The computation of these Nil-groups will be carried out in the next section.
\end{proof}

\section{Computing the Waldhausen Nil-groups.}

In this section we implement an approach suggested to us by F.T.
Farrell, and provide explicit computations for the Waldhausen
Nil-groups appearing in the previous section.  The results of this
section can be summarized in the following two theorems:

\begin{theorem}
For $i=0,1$, $NK_i(\mathbb ZD_3;C_1, C_2)\cong 0$ where $C_j=\mathbb
Z[D_6 \setminus D_3]$ is the $\mathbb ZD_3$-bimodule generated by
$D_6 \setminus D_3$ for $j=1,2$
\end{theorem}

\begin{theorem}
For $i=0,1$, $NK_i(\mathbb ZD_2;B_1, B_2)\cong \bigoplus _\infty
\mathbb Z/2$, where $B_j=\mathbb Z[D_2 \times Z/2 \setminus D_2]$ is
the $\mathbb ZD_2$ bi-module generated by $D_2 \times \mathbb Z/2
\setminus D_2$ for $j=1,2$
\end{theorem}

The proof of these theorems will be based on the following observations:

\begin{enumerate}
\item[(1)] There exists a {\it surjection}:

$$2\cdot NK_i(\mathbb Z[D_k]) \rightarrow NK_i(\mathbb ZD_k;M_1, M_2)$$
from the direct sum of two copies of the Bass Nil-groups
$NK_i(\mathbb Z[D_k])$ to the Waldhausen Nil-groups for $i=0,1$ ($k=2,3$).
$M_1,M_2$ are suitably defined bi-modules that specialize to
$B_1,B_2$ and $C_1,C_2$ when $k=2,3$ respectively.

\vspace{.2cm}
\item[(2)] The surjection above is an {\it injection} on each individual
copy of $NK_i(\mathbb Z[D_k])$.
\end{enumerate}

Let us now outline how Theorems 5.1 and 5.2 follow from these observations.  In the case where $k=3$, we know that
$NK_i(\mathbb Z[D_3])\cong 0$ ($i=0,1$), and the surjectivity of the map
immediately tells us that the Waldhausen Nil-groups $NK_i(\mathbb
ZD_3;C_1, C_2)$ vanish.  On the other hand, when $k=2$, we establish
that $NK_i(\mathbb Z[D_2])\cong \bigoplus _\infty \mathbb Z/2$
($i=0,1$), an infinite countable sum of $\mathbb Z/2$.  This is done
by establishing that $NK_i(\mathbb Z[D_2])\cong NK_{i+1}(\mathbb
F_2[\mathbb Z/2])$ ($i=0,1$), and appealing to computations of
Madsen \cite{Mad95} for the latter groups.  Once we have this result for
$NK_i(\mathbb Z[D_2])$, the surjectivity and injectivity statements
above combine to force $NK_i(\mathbb ZD_2;B_1, B_2)\cong \bigoplus
_\infty \mathbb Z/2$ ($i=0,1$).

%%%%%%%%%%%%%%%%%%%%%%%%%%%%%%%%%%%%%%%%%%%%%%%%%%%%%%%%%%

Let us start by recalling some well-known facts about infinity virtually cyclic groups.  If $\g$ is a virtually infinite cyclic group not of the form $G \rtimes \mathbb Z$, then $\g$ always maps epimorphically onto the infinite dihedral group $D_{\infty}$ with finite kernel (see \cite[Lemma 2.5]{FJ95}). Using this epimorphism,  Farrell and Jones in \cite{FJ95} constructed a stratified fiber bundle $\rho_ E: E \rightarrow X$,  where $E$ is a closed manifold with $\pi_1 E= \g$ and such that, for each point $x \in X$,
$\pi_1\rho^{-1}_E(x)=G_x$ or  $G_x \rtimes \mathbb Z$, where $G_x$ is a finite group.

\vskip 10pt

Next, we briefly describe the control space $X$ as constructed by Farrell and Jones in \cite[Section 2]{FJ95}.

Recall the model for the real hyperbolic $\mathbb H^2$ discussed in \cite[Section 2]{FJ86}, i.e.,
$\mathbb H^2 = \{(x,y) \in \mathbb R^2\,|\,-1 < y < 1\},$
with the Riemannian metric defined in \cite[pg. 23]{FJ95}. 

Let $S\mathbb H^2 \rightarrow \mathbb H^2$ be the tangent sphere bundle consisting of unit length vectors, and let   $S^+\mathbb H^2 \rightarrow \mathbb H^2$ be the subbundle of  $S\mathbb H^2 \rightarrow \mathbb H^2$ whose fiber over $x \in \mathbb H^2$ is $S^+_x \mathbb H^2$ ({\it the asymptotic northern hemisphere} defined in \cite[0.12]{FJ86}). The natural action of $D_{\infty}$ on $\mathbb H^2$ determined by $D_{\infty} \subseteq \iso(\mathbb H^2)$ induces an action of $D_{\infty}$ on $S^+\mathbb H^2$.

The orbit space $S^+\mathbb H^2/D_{\infty}$ has a stratification with six strata. The lowest strata are $H$, $V_0$, $V_1$, the intermediate strata are $B_0$, $B_1$ and $T$ is the top stratum. Using the canonical quotient map 
\[
p:S^+\mathbb H^2 \rightarrow S^+\mathbb H^2/D_{\infty},
\]
these strata are defined in \cite[pg. 24]{FJ95} .

Here the key point is to note that $H$ is diffeomorphic to the circle $S^1$ while $V_0$ and $V_1$ are both diffeomorphic to $\mathbb R$. Also $p(\partial S^+\mathbb H^2)$ is diffeomorphic to the cylinder $S^1 \times \mathbb R$, and $H$ cuts this manifold into two connected components $B_0$ and $B_1$, i.e.
\[
B_0 \coprod B_1 = p(\partial S^+\mathbb H^2) - H.
\]
Note that each pair $(B_i \cup H, H)$, $i= 0, 1$, is diffeomorphic to the pair $(S^1 \times [0, +\infty), S^1 \times 0)$. Finally, the top stratum is the complement in $S^+\mathbb H^2/D_{\infty}$ of the union $V_0 \cup V_1 \cup p(\partial S^+\mathbb H^2)$.

The control space $X$ is defined to be the quotient space of $S^+\mathbb H^2/D_{\infty}$ obtained by identifying the stratum $H$ to a single point $\ast$. Let
\[
\rho: S^+\mathbb H^2/D_{\infty} \rightarrow X 
\]
denote the canonical quotient map. Then the stratification of $S^+\mathbb H^2/D_{\infty}$ induces a stratification of $X$ whose six strata are the images of the strata in $S^+\mathbb H^2/D_{\infty}$ under $\rho$. Since
\[
\rho: S^+\mathbb H^2/D_{\infty} - H \rightarrow X - \ast
\]
is a homeomorphism, we identify via $\rho$ the strata in $S^+\mathbb H^2/D_{\infty}$ different from $H$ with the corresponding strata in $X$ different from $\ast$. 

This concludes the description of the control space $X$. For specific details and the construction of the stratified fiber bundle $\rho_E: E \rightarrow X$ mentioned earlier, we refer the reader to \cite[Section 2]{FJ95}.

Farrell and Jones also proved  in \cite[Theorem 2.6]{FJ95}  that the group homeomorphism
\[
\pi_{i}(\mathcal A): \pi_i \mathcal P_\ast(E ; \rho_E) \longrightarrow  \pi_i \mathcal P_\ast(E)
\]
is an epimorphism for every integer $i$. Here $\mathcal P_{\ast}(E)$ is the spectrum of stable topological pseudoisotopies on $E$, $\mathcal P_{\ast}(E, \rho_E)$ is the spectrum of those stable pseudoisotopies which are controlled over $X$ via $\rho_E$, and $\mathcal A$ is the `assembly' map.

\vskip 10pt

Recall that Quinn \cite[Theorem 8.7]{Qu82} constructed a spectral $E^n_{s,t}$, abutting to $ \pi_{s+t} \mathcal P_{\ast}(E ; \rho_E)$ with $E^2_{s,t} = H_s(X ; \pi_t \mathcal P_{\ast}(\rho_E))$. Here $\pi_q \mathcal P_{\ast}(\rho_E)$, $q \in \mathbb Z$, denotes the stratified system of abelian groups over $X$ where the group above $x \in X$ is $\pi_q(\mathcal P_{\ast}(\rho^{-1}_E(x))$. Note that by Anderson and Hsiang's result (see  \cite[Theorem 3]{AH} and Section 4), $\pi_i \mathcal P_{\ast +2}(\rho^{-1}_E(x))=K_i(\mathbb Z \pi_1(\rho^{-1}_E(x)))$ for $i \leq -1$, $\tilde{K}_0(\mathbb Z \pi_1(\rho^{-1}_E(x)))$, for $i=0$, and $Wh(\mathbb Z \pi_1(\rho^{-1}_E(x)))$ for $i=1$. 

It is important to emphasize that for each $x \in X$, the fundamental group $\pi_1 \rho^{-1}_E(x)$ of the fiber is either finite or a semidirect product $F \rtimes \mathbb Z$ where $F$ is a finite subgroup of $\g$  (see \cite[Remark 2.6.8]{FJ95}). 

These facts together with the immediate consequences of the construction of the control space $X$ and the stratified fiber bundle $\rho_E: E \to X$ are extremely  useful in analyzing Quinn's spectral sequence. This analysis is then used to complete the proof of the {\it surjectivity} part in Theorems 5.1 and 5.2.

\begin{proof}[Proof of Theorem 5.1]
For $\g=D_3 \times D_{\infty} \cong D_6 \ast_{D_3} D_6$,  first note that the fundamental group $F$ of a fiber of $\rho_E$, i.e. $\pi_1 \rho^{-1}_E(x)$ is either the group $D_3 \times \mathbb Z$ or some finite subgroup  $F$ of $\g$. It is a fact that $K_q( \mathbb Z F)$ (for $q \leq -1$), $\tilde{K}_0(\mathbb Z F)$, and $Wh(F)$ all vanish when $F$ is a finite subgroup of $\g$ except when $F=D_6$. This fact can be obtained by combining results from \cite{C80}, \cite{Re76}, \cite{Ro94}, \cite{Or04}, \cite{Pe98}, \cite{CuR81}, \cite{O89} (see the arguments in Section 4). It is also a fact that $\tilde{K}_0(\mathbb Z D_6)$,  $Wh(D_6)$ both vanish, and  $K_{-1}(\mathbb ZD_6)=\mathbb Z$ (see \cite[p. 350]{Or04}, \cite[p. 274]{Pe98}).  To compute the lower algebraic $K$-theory of the ring $\mathbb Z[D_3 \times \mathbb Z]=\mathbb Z[D_3][\mathbb Z]$, we use the Fundamental Theorem of algebraic $K$-theory (see \cite{Bas68}). The following are the results found: 

\[
\begin{aligned}
Wh(D_3 \times \mathbb Z) &\cong NK_1(\mathbb Z[D_3]) \oplus NK_1(\mathbb Z[D_3]),\\
\tilde{K}_0(\mathbb Z[D_3 \times \mathbb Z]) &\cong  NK_0(\mathbb Z[D_3])  \oplus NK_0(\mathbb Z[D_3]),\\
K_q(\mathbb Z[D_3 \times \mathbb Z]) &\cong 0 \quad \text{for} \quad  q \leq -1,
\end{aligned}
\]

\noindent where the Nil-groups appearing in these computations are the Bass's Nil-groups.

Next, note that given the simple nature of the control space $X$, the spectral  sequence collapses at $E^2$ and for all $s, t \in \mathbb Z$
\[
E^2_{s,t} = H_s(X ; \pi_t \mathcal P_{\ast}(\rho_E)) \Longrightarrow \pi_{s+t} \mathcal P_{\ast}(E ; \rho_E) 
\]
Combining the above with the fact that $Wh_q(D_3)$ vanish for all $q \leq 1$, we have that for all $s, t \in \mathbb Z$
\[
\begin{aligned}
 H_s(X ; \pi_t \mathcal P_{\ast}(\rho_E))
              &=H_s(\ast; Wh_t(D_3\times \mathbb Z)) \oplus H_s(\mathbb R; Wh_t(D_6)) \oplus H_s(\mathbb R; Wh_t(D_6))\\
	&=Wh_t(D_3\times \mathbb Z) \oplus Wh_t(D_6) \oplus Wh_t(D_6).
\end{aligned}
\]
Consequently,
\[
\begin{aligned}
\pi_{-1} \mathcal P_{\ast}(E ; \rho_E) &\cong  NK_1(\mathbb Z[D_3])  \oplus NK_1(\mathbb Z[D_3]),\\
\pi_{-2} \mathcal P_{\ast}(E ; \rho_E) &\cong  NK_0(\mathbb Z[D_3])  \oplus NK_0(\mathbb Z[D_3]),\\
\pi_{-3} \mathcal P_{\ast}(E ; \rho_E) &\cong K_{-1}(\mathbb Z[D_3\times \mathbb Z]) \oplus K_{-1}(\mathbb Z[D_6]) \oplus K_{-1}( \mathbb Z[D_6]) \cong \mathbb Z \oplus \mathbb Z,\\
\pi_q \mathcal P_{\ast}(E ; \rho_E) &=0  \quad \text{when} \quad q \leq -4,
\end{aligned}
\]
On the other hand using the results found in Section 4 for the lower algebraic $K$-theory of the integral group ring of $D_3 \times D_{\infty}$, we have that 
\[
\begin{aligned}
\pi_{-1} \mathcal P_{\ast}(E) &\cong Wh(D_3 \times  D_{\infty}) \cong NK_1(\mathbb ZD_3; C_1, C_2)\\
\pi_{-2} \mathcal P_{\ast}(E) &\cong \tilde{K}_0(\mathbb Z[D_3 \times   D_{\infty}]) \cong  NK_0(\mathbb ZD_3 ; C_1, C_2)\\
\pi_{-3} \mathcal P_{\ast}(E) &\cong K_{-1}(\mathbb Z[D_3 \times D_{\infty}]) \cong \mathbb Z \oplus \mathbb Z\\
\pi_q \mathcal P_{\ast}(E) &\cong K_q(\mathbb Z[D_3 \times   D_{\infty}]) \cong 0 \quad \text{for} \quad  q \leq 4,
\end{aligned}
\]

\noindent where $C_i=\mathbb Z[D_6 \setminus D_3]$ is the
$\mathbb ZD_3$-bimodule generated by $D_6 \setminus D_3$ for
$i=1,2$. Here the Nil-groups appearing in these computations are the Waldhausen's Nil-groups.

Combining the above with \cite[Theorem 2.6]{FJ95} we obtain for $i=0, 1$ the desired  epimorphism
\[
2 \cdot NK_i(\mathbb Z[D_3]) \longrightarrow NK_i(\mathbb ZD_3, C_0, C_1) \longrightarrow 0
\]

\noindent Since $NK_i(\mathbb Z[D_3])=0$ for $i=0,1$ (see \cite{Ha87}), then it follows that for $i=0, 1$, 
$NK_i(\mathbb ZD_3, C_0, C_1)\cong 0$, completing the proof of Theorem 5.1. 
\end{proof}

Before proving Theorem 5.2, we prove the following Lemmas.

\begin{lemma}
\[
NK_0(\mathbb Z[D_2])\cong NK_{1}(\mathbb F_2[\mathbb Z/2]) \cong \bigoplus_{\infty} \mathbb Z/2.
\] 
\end{lemma}

\begin{proof}
Write $\mathbb Z[D_2]=\mathbb Z[\mathbb Z/2 \times \mathbb Z/2]$ as $\mathbb Z[\mathbb Z/2][\mathbb Z/2]$, and  consider the following Cartesian square
\[
\xymatrix@C=20pt@R=30pt{
\mathbb Z[\mathbb Z/2][\mathbb Z/2] \ar[d] \ar[r] & \mathbb Z[\mathbb Z/2]
     \ar[d] \\
\mathbb Z[\mathbb Z/2] \ar[r] & \mathbb F_2[\mathbb Z/2]}
\]      

\noindent This Cartesian square yields the Mayer-Vietories sequence 

\begin{equation*}
\begin{split}
\cdots \rightarrow NK_1(\mathbb Z[\mathbb Z/2]) \oplus NK_1(\mathbb Z[\mathbb Z/2]) &\rightarrow NK_1(\mathbb F_2[\mathbb Z/2]) \rightarrow NK_0(\mathbb Z[\mathbb Z/2][\mathbb Z/2]) \rightarrow \\
&\rightarrow NK_0(\mathbb Z[\mathbb Z/2]) \oplus NK_0(\mathbb Z[\mathbb Z/2]) \rightarrow \cdots
\end{split}
\end{equation*}

\noindent Since $NK_i(\mathbb Z[\mathbb Z/2])=0$, for $i=0, 1$ (see \cite{Ha87}), we obtain that 
$NK_0(\mathbb Z[D_2])\cong NK_{1}(\mathbb F_2[\mathbb Z/2])$. In \cite{Mad95}, Madsen showed 
that $NK_{1}(\mathbb F_2[\mathbb Z/2])\cong \bigoplus_{\infty} \mathbb Z/2$, giving us the 
desired result.

\end{proof}

\begin{lemma}
\[
NK_1(\mathbb Z[D_2]) \cong NK_{2}(\mathbb F_2[\mathbb Z/2]) \cong \bigoplus_{\infty} \mathbb Z/2.
\] 
\end{lemma}

\begin{proof}
As before, consider the following Cartesian square
\[
\xymatrix@C=20pt@R=30pt{
\mathbb Z[\mathbb Z/2][\mathbb Z/2] \ar[d] \ar[r] & \mathbb Z[\mathbb Z/2]
     \ar[d] \\
\mathbb Z[\mathbb Z/2] \ar[r] & \mathbb F_2[\mathbb Z/2]}
\]      
\noindent As mentioned in Lemma 5.3, since $NK_1(\mathbb Z[\mathbb Z/2])$ vanishes,  the Mayer-Vietories sequence 
yields the epimorphism
\[
NK_2(\mathbb F_2[\mathbb Z/2]) \rightarrow NK_1(\mathbb Z[\mathbb Z/2][\mathbb Z/2])   \longrightarrow 0 
\]

\noindent Now,  in \cite{Mad95}, Madsen has shown that $NK_{2}(\mathbb F_2[\mathbb Z/2])\cong 
\bigoplus_{\infty} \mathbb Z/2$.  Combining this result with the well-known fact that the Nil groups are either  
trivial or infinitely generated \cite{F77}, we conclude that up to isomorphism, either
\[
NK_1(\mathbb Z[D_2]) \cong\bigoplus_{\infty} \mathbb Z/2
\]
\noindent or the group is trivial.  But the group $NK_1(\mathbb Z[\mathbb Z/2][\mathbb Z/2])$ is known
to be non-trivial (see \cite{Bas68}), completing the proof of the Lemma.
\end{proof}

\begin{proof}[Proof of Theorem 5.2]
Let $\g=D_2 \times D_{\infty} \cong D_2 \times \mathbb Z/2 \ast_{D_2} D_2 \times \mathbb Z/2$. In this case  
note that the fundamental group $F$ of a fiber of $\rho_E$, i.e. $\pi_1 \rho^{-1}_E(x)$ is either the group 
$D_2 \times \mathbb Z$ or some finite subgroup  $F$ of $\g$. It is a fact that $K_i( \mathbb Z F)$ 
(for $i \leq -1$), $\tilde{K}_0(\mathbb Z F)$, and $Wh(F)$ all vanish when $F$ is a finite subgroup of $\g$. 
This fact can be obtained by combining results from \cite{C80}, \cite{Re76}, \cite{Ro94}, \cite{Or04}, 
\cite{Pe98}, \cite{CuR81}, \cite{O89} (see Section 4).  To compute the lower algebraic $K$-theory of the ring 
$\mathbb Z[D_2 \times \mathbb Z] = \mathbb Z[D_2][\mathbb Z]$,  we apply the Fundamental Theorem of 
algebraic $K$-theory to obtain:
\[
\begin{aligned}
Wh(D_2 \times \mathbb Z) &\cong NK_1(\mathbb Z[D_2]) \oplus NK_1(\mathbb Z[D_2])\\
\tilde{K}_0(\mathbb Z[D_2 \times \mathbb Z]) &\cong  NK_0(\mathbb Z[D_2])  \oplus NK_0(\mathbb Z[D_2])\\
K_q(\mathbb Z[D_2 \times \mathbb Z]) &\cong 0 \quad \text{for} \quad  q \leq -1,
\end{aligned}
\]
where the Nil-groups appearing in these computations are the Bass's Nil-groups.

Once again, given the simple nature of the control space $X$ in this case, the spectral sequence collapses at $E^2$ and for all $s, t \in \mathbb Z$
\[
E^2_{s,t} = H_s(X ; \pi_t \mathcal P_{\ast}(\rho_E)) \Longrightarrow \pi_{s+t} \mathcal P_{\ast}(E ; \rho_E) 
\]
Combining the above with the fact that $Wh_{q}(D_2)$ and $Wh_q(D_2 \times \mathbb Z/2)$ are both trivial for all $q \leq 1$, we have that for all $s, t \in \mathbb Z$
\[
 H_s(X ; \pi_t \mathcal P_{\ast}(\rho_E)) =Wh_t(D_2\times \mathbb Z) 
\]
Consequently, we obtain:
\[
\begin{aligned}
\pi_{-1} \mathcal P_{\ast}(E ; \rho_E) &\cong  NK_1(\mathbb Z[D_2])  \oplus NK_1(\mathbb Z[D_2]),\\
\pi_{-2} \mathcal P_{\ast}(E ; \rho_E) &\cong  NK_0(\mathbb Z[D_2])  \oplus NK_0(\mathbb Z[D_2]),\\
\pi_q \mathcal P_{\ast}(E ; \rho_E) &=0  \quad \text{when} \quad q \leq -3.
\end{aligned}
\]
On the other hand, using the results found in Section 4 for the lower algebraic $K$-theory of $D_2 \times D_{\infty}$, we have that
\[
\begin{aligned}
\pi_{-1} \mathcal P_{\ast}(E) &\cong Wh(D_2 \times  D_{\infty}) \cong NK_1(\mathbb ZD_2; B_1, B_2)\\
\pi_{-2} \mathcal P_{\ast}(E) &\cong \tilde{K}_0(\mathbb Z[D_2 \times   D_{\infty}]) \cong  NK_0(\mathbb ZD_2 ; B_1, B_2)\\
\pi_{q} \mathcal P_{\ast}(E) &\cong K_i(\mathbb Z[D_2 \times   D_{\infty}]) \cong 0 \quad \text{for} \quad  q \leq -3,
\end{aligned}
\]
where $B_i=\mathbb Z[D_2 \times Z/2 \setminus D_2]$ is the $\mathbb ZD_2$ bi-module generated by $D_2
\times \mathbb Z/2 \setminus D_2$ for $i=1,2$. Recall that the Nil-groups appearing in these computations are the Waldhausen's Nil-groups.

Combining the above with \cite[Theorem 2.6]{FJ95} we obtain for $i=0, 1$ the  epimorphism
\[
2 \cdot NK_i(\mathbb Z[D_2]) \longrightarrow NK_i(\mathbb ZD_2, B_1, B_2) \longrightarrow 0
\]

%%%%%%%%%%%%%%%%%%%%%%%%%%%%%%%%%%%%%%%%%%%%%%%%%%%%%%%
%%%%%%%%%%%%%%%%%%%%%%%%%%%%%%%%%%%%%%%%%%%%%%%%%%%%%%%%
Next, recall that a monomorphism of group $\sigma: \Pi \rightarrow \g$ induces transfer maps $\sigma^{\ast}: Wh(\g) \rightarrow  Wh(\Pi)$ and $\sigma^{\ast}: \tilde K_0( \mathbb Z\g) \rightarrow  \tilde K_0(\mathbb Z\Pi)$ provided the image of $\sigma$ have finite index in $\g$ (see \cite{FH78}).  

Let $\Pi = G \times \mathbb Z$, where  $G=D_2$. By the above paragraph, the monomorphism $\sigma$ in the exact sequence
\[
0 \longrightarrow \Pi \overset{\sigma}{ \longrightarrow} \g \longrightarrow \mathbb Z/2 \longrightarrow 0
\]
 induces the transfer maps:
\[
\sigma^{\ast}:  Wh(G \times D_{\infty}) \longrightarrow Wh(G \times \mathbb Z)
\]

\[
\sigma^{\ast}: \tilde K_0(\mathbb Z[G \times D_{\infty}]) \longrightarrow \tilde K_0(\mathbb Z[D_2 \times \mathbb Z])
\]
Let $\tau$  be the standard involution of $G$ given by $\tau(g)=g^{-1}$ for $g \in G$. This involution can be extended to an involution of $\mathbb Z[G \times \mathbb Z]=\mathbb Z[G][\mathbb Z]$ (see \cite[pg. 209]{FH68}).  By abuse of 
notation, we will also denote this new involution by $\tau$. 

In the direct sum decomposition given by the {\it Fundamental Theorem of algebraic $K$-theory} of the algebraic $K$-groups of the ring $\mathbb Z[G][\mathbb Z]$:

\[
K_i(\mathbb Z[G][\mathbb Z]) \cong K_{i-1} (\mathbb Z[G]) \oplus K_i(\mathbb Z[G]) \oplus NK_i(\mathbb Z[G]) \oplus NK_i(\mathbb Z[G]),
\]
$\tau_{\ast}$ interchanges the two direct summands $NK_i(\mathbb Z[G])$ (see \cite[Proposition 20]{FH68}).
Now recall that  $(\sigma^{\ast} \circ \sigma_{\ast})(x) = x +\tau_{\ast}(x)$ for all $x \in Wh_i(D_2 \times \mathbb Z)$, $i=0,1$; we will denote this map by $\psi$.  This yields the following commutative diagram for $i=0, 1$:
\[
\xymatrix{
 Wh_i(D_2 \times \mathbb Z) \ar[dr]_{\psi}  \ar[r]^{\sigma_{\ast}} & Wh_i(D_2 \times D_{\infty})
     \ar[d]^{\sigma^{\ast}} \\ 
&Wh_i(D_2 \times \mathbb Z)}
\]      
Using our earlier computations, we obtain for $i=0, 1$ the commutative diagram:
\[
\xymatrix{
 NK_i(\mathbb Z[D_2]) \oplus NK_i(\mathbb Z[D_2])   \ar[dr]_{\psi}  \ar[r]^{\sigma_{\ast}} & NK_i(\mathbb ZD_2, B_1, B_2)
     \ar[d]^{\sigma^{\ast}} \\ 
&NK_i(\mathbb Z[D_2]) \oplus NK_i(\mathbb Z[D_2])}
\]      
We now define the map $\varphi$ by:
\[
\varphi:  NK_i(\mathbb Z[D_2]) \xrightarrow{\text{(\text{Id}, 0)}} NK_i(\mathbb Z[D_2]) \oplus  NK_i(\mathbb Z[D_2])
\]

Note that $\sigma^{\ast} \circ \sigma_{\ast}$ is injective on  $\varphi(NK_i(\mathbb Z[D_2]))$.  It follows that $\sigma_{\ast}$ is injective on the first summand of the direct sum $NK_i(\mathbb Z[D_2]) \oplus NK_i(\mathbb Z[D_2])$. 
Note that an identical argument can be used to show that $\sigma^{\ast} \circ \sigma_{\ast}$ is injective on the
second summand (though it is easy to see that the image of the first and the second summand coincide).

\vspace{.3cm}
For $i=0,1$, this gives us the desired monomorphism:
\[
0 \longrightarrow NK_i(\mathbb Z[D_2]) \longrightarrow NK_i(\mathbb ZD_2, B_1, B_2).
\]
Combining this monomorphism with the epimorphism found earlier,
and applying Lemma 5.3 and  Lemma 5.4, it follows that $i=0,1$
\[
NK_i(\mathbb ZD_2, B_1, B_2) \cong \bigoplus_{\infty} \mathbb Z/2,
\]
This completes the proof of Theorem 5.2.
\end{proof}

%%%%%%%%%%%%%%%%%%%%%%%%%%%%%%%%%%%%%%%%%%%%%%%%%%%%%%%%%%%%

\end{document}